\documentclass[12pt]{article}
\usepackage{amsfonts}
\usepackage{amssymb}
\usepackage{amscd, color}
\usepackage{hyperref}
\usepackage{amsmath}
\newcommand{\BEQ}{\begin{equation}}
\newcommand{\EEQ}{\end{equation}}
\newtheorem{Th}{Theorem}

\newtheorem{Lem}{Lemma}
\newtheorem{Prop}{Proposition}
\newtheorem{Cor}{Corollary}

\newtheorem{Ex}{Example}
\newtheorem{Def}{Definition}
\newtheorem{Rem}{Remark}
\def\nn{\nonumber}
\def\one#1{#1^{\raise5pt\hbox{$\scriptstyle\!\!\!\!1$}}\,{}}
\def\two#1{#1^{\raise5pt\hbox{$\scriptstyle\!\!\!\!2$}}\,{}}

\def\bea{\begin{eqnarray}}
\def\eea{\end{eqnarray}}

\newcommand{\CC}{{{\mathbb{C}}}}

\def\l{\lambda}

\def\T{\mathcal{T}}

\def\p{\partial_z}

\def\tpsi{\tilde \psi}

\def\GT{Manin}
\def\MMs{Manin\ matrices}
\def\MM{Manin\ matrix}

\def\BX{$\Box$}
\def\PRF{ {\bf Proof.} }
\def\SPRF{ {\bf Sketch of the Proof.} }
\def\ms{ n } 
\def\mpn{ k } 
\def\GIOPS{ Gurevich, Isaev, Ogievetsky, Pyatov, Saponov papers \cite{GIOPS95}-\cite{GIOPS05}}
\def\OkOlsh{\href{http://lanl.arxiv.org/abs/q-alg/9605042}{[OkounkovOlshansky96]}} 

\def\GR97{\href{http://arxiv.org/abs/q-alg/9705026}{[GR97]}}
\def\KL94{\href{http://arxiv.org/abs/q-alg/9411194}{[KL94]}}

\def\MNO{\href{http://arxiv.org/abs/math.QA/0211288}{[MNO94]}}
\def\Kir{\href{http://www.ams.org/distribution/mmj/vol1-1-2001/kirillov.pdf}{[Kir00]}}
\def\GLZ{\href{http://arxiv.org/abs/math.QA/0303319}{[GLZ03]}}
\def\GIOPSSupB{\href{http://arxiv.org/abs/math.QA/0508506}{[GIOPS05]}}
\def\GIOPSW{\href{http://lanl.arxiv.org/abs/math/9809170}{[GIOPS98]}}
\def\MTVCap{\href{http://arxiv.org/abs/math.QA/0610799}{[MTV06]}}
\def\Ok961{\href{http://arxiv.org/abs/q-alg/9602028}{[Ok96]}}
\def\CTBig{\href{http://arxiv.org/abs/hep-th/0604128}{[CT06]}}
\def\T04{\href{http://arxiv.org/abs/hep-th/0404153}{[Ta04]}}
\def\MTVShap{\href{http://arxiv.org/abs/math.AG/0512299}{[MTV05]}}
\def\KonvalinkaA{\href{http://arxiv.org/abs/math.CO/0703213}{[Ko07B]}}
\def\KonvalinkaB{\href{http://arxiv.org/abs/math.CO/0703203}{[Ko07A]}}
\def\SklA{\href{http://arxiv.org/abs/hep-th/9211111}{[Sk92]}}
\def\SklB{\href{http://arxiv.org/abs/solv-int/9504001}{[Sk95]}}
\def\FMA{\href{http://arxiv.org/abs/nlin.SI/0306005}{[FM03A]}}
\def\FMB{\href{http://arxiv.org/abs/nlin.SI/0306008}{[FM03B]}}

\def\RybA{\href{http://arxiv.org/abs/math.RT/0606380}{[Ry1]}}
\def\RybB{\href{http://arxiv.org/abs/math.QA/0608588}{[Ry3]}}
\def\CRT{\href{http://arxiv.org/abs/hep-th/0404106}{[CRyT04]}}
\def\BT07{\href{http://arxiv.org/abs/hep-th/0703124}{[BT07]}}
\def\DubrovinM2{\href{http://arxiv.org/abs/math/0311261}{[DM03]}}
\def\Molev02{\href{http://arxiv.org/abs/math.QA/0211288}{[Mo02]}}
\def\CTKZ{\href{http://arxiv.org/abs/hep-th/0607250}{[CT06b]}}
\def\EnriquezRubtsov01{\href{http://arxiv.org/abs/math.AG/0112276}{[ER01]}}
\def\BabelonTalon02{\href{http://arxiv.org/abs/hep-th/0209071}{[BT02]}}
\def\Nazarov91{\href{http://www.springerlink.com/content/rt5r3313732p48j7}{[Na91]}}
\def\GelfandGelfandRetakhWilson02{\href{http://arxiv.org/abs/math.QA/0208146}{[GGRW02]}}
\def\CTGoper{\href{http://arxiv.org/abs/hep-th/0409007}{[CT04]}}

\def\EtingofPak06{\href{http://arxiv.org/abs/math/0608005}{[EP06]}}

\def\HaiLorenz06{\href{http://arxiv.org/abs/math/0603169}{[HL06]}}

\def\KonvalinkaPak06{\href{http://arxiv.org/abs/math/0607737}{[KPak06]}}

\def\FoataHan06{\href{http://arxiv.org/abs/math/0603463}{[FH06]}}

\def\Sklyanin92{\href{http://arxiv.org/abs/hep-th/9212076}{[Sk92b]}}

\def\AHH96{\href{http://arxiv.org/abs/solv-int/9603001}{[AHH96]}}

\def\CTHA{\href{http://arxiv.org/abs/hep-th/0303069}{[CT03A]}}

\def\CTHB{\href{http://arxiv.org/abs/hep-th/0309059}{[CT03B]}}

\def\Sah06{\href{http://arxiv.org/abs/math-ph/0609067}{[Sa06]}}

\def\GGM97{\href{http://arxiv.org/abs/hep-th/9707120}{[GGM97]}}

\def\NO{\href{http://arxiv.org/abs/q-alg/9507003}{[NO95]}}

\def\BartocciFalquiPedroni03{\href{http://lanl.arxiv.org/abs/nlin/0307021}{[BaFP03]}}

\def\PetreraRagnisco07{\href{http://lanl.arxiv.org/abs/math-ph/0703044}{[PR07]}}

\def\AdlerMoerbeke80{\href{http://www.sciencedirect.com}{[AvM80]}}

\def\ReymanSemenov79{\href{http://www.springerlink.com/content/n7360721350h7621}{[RS79]}}

\def\Beaville90{\href{http://www.springerlink.com/content/jn12352983882un4}{[Be90]}}

\def\AdamsHarnadPreviato88{\href{http://projecteuclid.org/DPubS/Repository/1.0/Disseminate?view=body&id=pdf_1&handle=euclid.cmp/1104161743}{[AHP88]}}

\def\CFR{\href{http://lanl.arxiv.org/abs/0710.4971}{[CFRy07]}}

\def\Frenkel95{\href{http://lanl.arxiv.org/abs/q-alg/9506003}{[Fr95]}}

\def\Reshetikhin92{\href{http://www.springerlink.com/content/t718750378552740}{[Re92]}}

\def\ManinA{\href{http://www.numdam.org/item?id=AIF_1987__37_4_191_0}{[Man87]}}

\def\FFR{\href{http://arxiv.org/abs/hep-th/9402022}{[FFR94]}}

\def\UmedaW{\href{http://www.springerlink.com/content/2l45hauppwnjc6be/}{[U03]}}

\def\UmedaN{\href{http://www.ams.org/proc/1998-126-11/S0002-9939-98-04557-2/home.html}{[U98]}}

\def\ItoN{\href{http://www.ingentaconnect.com/content/els/00218693/1998/00000208/00000002}{[Ito98]}}

\def\BT{\href{http://arxiv.org/abs/hep-th/0703124}{[BTa07]}}

\def\OrtizSommaDukelskyRombouts04{\href{http://arxiv.org/abs/cond-mat/0407429}{[OSDR04]}}

\def\GNR{\href{http://arxiv.org/abs/hep-th/9901089}{[GNR99]}}

\def\FaddeevTakhdadjanSklyanin79{\href{http://www.springerlink.com/content/uk1268x006041782/?p=6121269a4db64b09bba3c959e06b3a30&pi=1}
{[FTS79]}}

\def\AdamsHarnadHurtubise93{\href{http://projecteuclid.org/DPubS?service=UI&version=1.0&verb=Display&handle=euclid.cmp/1104253285}
{[AHH93]}}

\def\KZ{\href{http://dx.doi.org/10.1016/0550-3213(84)90374-2}{[KZ84]}}

\def\KulishSklyanin81{\href{http://www.springerlink.com/content/750u2k2xg12w8886/}{[KS81]}}

\def\ArnaudonAvanCrampeDoikouFrappatRagoucy04{\href{http://arxiv.org/abs/math-ph/0406021}
{[AACDFR04]}}

\def\Gekhtman{\href{http://projecteuclid.org/DPubS?service=UI&version=1.0&verb=Display&handle=euclid.cmp/1104272160}
{[Ge95]}}

\def\Hayashi88{\href{http://www.springerlink.com/content/wp5nwm73g35147u5/}{[Ha88]}}

\def\Kontsevich97{\href{http://arxiv.org/abs/q-alg/9709040}{[K97]}}

\def\Nekrasov95{\href{http://lanl.arxiv.org/abs/hep-th/9503157}{[Ne95]}}
\def\EnriquezRubtsovA95{\href{http://lanl.arxiv.org/abs/alg-geom/9503010}{[ER95]}}

\def\OkounkovB96{\href{http://lanl.arxiv.org/abs/q-alg/9602027}{[Ok96B]}}

\def\NazarovB96{\href{http://lanl.arxiv.org/abs/q-alg/9601027}{[Na96]}}

\def\Hitchin87{\href{http://projecteuclid.org/euclid.dmj/1077305506}{[Hi87]}}

\def\RRT02{\href{http://www.springerlink.com/content/m6842k7r0u0q461u}{[RRT02]}}

\def\RRTb05{\href{http://dx.doi.org/10.1016/j.jalgebra.2005.01.002}{[RRT05]}}

\def\LT07{\href{http://dx.doi.org/10.1016/j.jpaa.2006.03.017}{[LT07]}}

\begin{document}
\hfill ITEP-TH-45/07

\vskip 1.5cm

\centerline{\Large Manin matrices and Talalaev's formula.}

\vskip 1.5cm
 \centerline{A. Chervov \footnote{E-mail:
chervov@itep.ru}
G. Falqui  \footnote{E-mail:gregorio.falqui@unimib.it}
}
\centerline{\sf Institute for Theoretical and Experimental Physics}

\centerline{\sf Universit\`a degli Studi di Milano - Bicocca }

\vskip 0.5cm
\centerline{\large \bf Abstract}
\vskip 1cm
In this
paper we study properties of Lax   and Transfer matrices associated
with quantum integrable systems. Our point of view stems from the
fact that their elements satisfy special commutation properties,
considered by Yu. I. Manin some twenty years ago at the beginning of
Quantum Group Theory. They are the commutation properties of matrix
elements of linear homomorphisms between polynomial rings; more
explicitly they read: 1) elements in the same column commute; 2)
commutators of the cross terms are equal: $[M_{ij}, M_{kl}]=[M_{kj},
M_{il}]$ (e.g. $[M_{11}, M_{22}]=[M_{21}, M_{12}]$).

The main aim of this paper is twofold: on the one hand we observe
and prove that such matrices (which we call {\em
  Manin} matrices for short) behave
almost as well as matrices with commutative elements. Namely
theorems of linear algebra (e.g., a natural definition of the
determinant, the Cayley-Hamilton theorem, the Newton identities and
so on and so forth) have a straightforward counterpart in the case
of Manin matrices.

On the other hand, we remark that such matrices are somewhat
ubiquitous in the theory of quantum integrability. For instance,
Manin matrices (and their q-analogs) include matrices satisfying the
Yang-Baxter relation "RTT=TTR" and the so--called Cartier-Foata
matrices. Also, they enter
\href{http://arxiv.org/abs/hep-th/0404153}
{Talalaev's} remarkable formulas: $det(\partial_z-L_{Gaudin}(z))$,
$det(1-e^{-\p}T_{Yangian}(z))$ for the "quantum spectral curve",
and appear in separation of variables problem and Capelli
identities.

We show that theorems of linear algebra, after being established for
such matrices, have various applications to quantum integrable
systems and Lie algebras, e.g in the construction of
new generators in $Z(U_{crit}(\widehat{gl_n}))$ (and, in general, in the
construction of quantum
conservation laws), in the Knizhnik-Zamolodchikov equation,
and in the problem of Wick ordering.

We also discuss applications to the
separation of variables problem, new Capelli identities and the
Langlands correspondence.

\newpage
\tableofcontents

\section{Introduction and Summary of the Results}
It is well-known that  matrices with generically noncommutative
elements do not admit a natural construction of the determinant and
basic theorems of the linear algebra fail to hold true. On the other
hand, matrices with noncommutative entries play a basic role in the
theory of quantum integrability
(see, e.g., \FaddeevTakhdadjanSklyanin79), 
in Manin's theory of "noncommutative symmetries" \cite{Manin}, and
so on and so forth. It is fair to say that recently D. Talalaev
{\T04} made a kind of breakthrough in quantum integrability,
defining the "quantum characteristic polynomial" or "quantum
spectral curve" $det(\p-L(z))$ for Lax matrices satisfying rational 
$R$-matrix commutation relations (e.g., Gaudin systems)\footnote{See
\CTBig, \MTVCap, \MTVShap\ for some applications to the Bethe
ansatz, separation of variables, the Langlands correspondence, the
Capelli identities, real algebraic geometry and related topics.}.

The first and basic observation  of the present paper is the
following: the {\em quantum Lax  matrix} $(\p-L_{gl_n-Gaudin}(z))$
entering\footnote{Remark that $L_{gl_n-Gaudin}(z)$ is not a \MM; to
insert $\p$ is an insight due to D.Talalaev; the operator
$(\p-L_{gl_n-Gaudin}(z))$ is also related to the
Knizhnik-Zamolodchikov equation
see section \ref{KZ}
} Talalaev's formula, as well as a suitable modification of the
transfer matrix of the (XXX) Heisenberg chain, matches the simplest
case of Manin's considerations (i.e. are \MMs~ in the sense
specified by the Definition \ref{D1} below). Further we prove that
many {\em results of commutative linear algebra} can be applied with
minor modifications in the case of "\MMs" and derive applications.

We will consider the simplest case of those considered by Manin,
namely - in the present paper - we will restrict ourselves to the
case of commutators, and not of (super)-$q$-commutators, etc.
Let us mention that \MMs~ are defined, roughly speaking, by 
{\bf half} of the relations
of the corresponding quantum group $Fun_q(GL(n))$ and taking $q=1$ (see remark \ref{QGrem} page
\pageref{QGrem}).
\begin{Def}\label{D1}
Let $M$ be a $n\times n'$ matrix with elements $M_{ij}$ in
a noncommutative ring $\mathcal{R}$. We will call $M$ a \MM\ if the
following two conditions hold:
\begin{enumerate}
\item Elements in the same column commute between themselves.
\item Commutators of cross terms of $2\times 2$ submatrices of $M$ are equal:
\begin{equation}\label{Cross-term-propert} [M_{ij}, M_{kl}]=[M_{kj}, M_{il}],\>
\forall\, i,j,k,l \text{ e. g. } [M_{11}, M_{22}]=[M_{21}, M_{12}].
\end{equation}
\end{enumerate}
\end{Def}
A more intrinsic definition of \MMs\ via coactions on polynomial and
Grassmann algebras will be recalled in Proposition \ref{proposizione} page \pageref{proposizione}
below. (Roughly speaking variables $\tilde x_i = \sum_j M_{ij}x_j$ commute
among themselves if and only if $M$ is a \MM, where
$x_j$ are commuting variables, also commuting with elements of $M$).

As we shall see in Section \ref{sect3}, the following properties hold:
\begin{itemize}
\item
$\p - L_{gl_n-Gaudin}(z)$ is a \MM. Here $L_{gl_n-Gaudin}(z)$ is the
Lax matrix for the Lie algebra $gl_n[t]$, i.e. a matrix
valued-generating function for generators of $gl_n[t]$, (in physics
$L_{gl_n-Gaudin}(z)$ is the Lax matrix for the Gaudin integrable
system). (Section \ref{Gaud-Manin-ss}).
\item
$e^{-\p} T_{gl_n-Yangian}(z)$ is a \MM. Here $T_{gl_n-Yangian}(z)$
is the Lax (or "transfer") matrix for the Yangian algebra $Y(gl_n)$,
i.e. a matrix valued-generating function for the generators of
$Y(gl_n)$, (in physics $T_{gl_n-Yangian}(z)$ is the Lax matrix for
the integrable systems with the Yangian symmetry: the Heisenberg XXX
spin chain, the Toda system, etc.)  (Section \ref{Yang-ss}).
\end{itemize}
Furthermore, we shall see that \MMs\ enter other topics in quantum
integrability, e.g, the quantum separation of variables theory
(section \ref{SV-s}) and Capelli identities (section
\ref{Cap-sssect}).
\\
Among the basic statements of linear algebra we establish (in a form
{\em identical} to the commutative case) for \MMs, let us mention:
\begin{itemize}
\item
Section \ref{Inv-c-ss}:  The inverse of a \MM~ $M$
is again Manin;
furthermore $det(M^{-1})=(det~M)^{-1}$;
\item Section \ref{Schur-ss}: the formula for the determinant of block matrices:
\[
 det\left(\begin{array}{cc}
A & B\\
C & D \\
\end{array}\right)=det(A)det(D-CA^{-1}B)=det(D)det(A-BD^{-1}C);\]
Also, we show that  $ D-CA^{-1}B$, $A-BD^{-1}C$ are \MMs;
\item
Section \ref{CH-ss}: the Cayley-Hamilton theorem: $det(t-M)|_{t=M}=0$;
\item
section \ref{Newton-ss}: the Newton identities between $TrM^k$ and
coefficients of $det(t+M)$).
\end{itemize}
This extends some results previously obtained in the literature:
Manin has defined  the determinant,
proved Cramer's inversion rule, Laplace formulas, as well as Plucker
identities.
In \GLZ ~
the MacMahon-Wronski formula
was proved; {\KonvalinkaA,\KonvalinkaB\ }
contains Sylvester's identity and the Jacobi ratio's theorem, along
with partial results on an inverse matrix and block
matrices.\footnote{These authors actually considered more general classes of matrices
}
One of our point is to
present applications to the theory of Quantum integrable systems.
\\
\paragraph{Applications to integrability and Lie theory.} We derive
various applications of the above-mentioned linear algebra statements
to quantum integrable systems, to the Knizhnik-Zamolodchikov
equation, as well as Yangians and Lie algebras. Namely:
\begin{itemize}
\item Section \ref{KZ}: a new proof  of the correspondence (\CTGoper)\
between solutions of the  Knizhnik-Zamolodchikov equation and
solutions $Q(z)$ of $det(\p-L(z)) Q(z)=0$. The last equation plays a
key role in the solution of the Gaudin type integrable systems. In
the $sl_2$ case it is known in the literature as "Baxter's T-Q
equation".
\item Section \ref{ERBT-ss}: the observation that
the inverse of a \MM~ is again a \MM\ yields, in a very simple way,
some recent result on a general construction about quantum
integrable systems (\EnriquezRubtsov01,\BabelonTalon02), that, in a
nutshell, says that (quantum)separability implies commutativity of the quantum
Hamiltonians.
\item Section \ref{Cap-sssect}: a new proof
of the generalized Capelli identities \MTVCap\
for the  Lie algebra $gl[t]$ (or "the Gaudin integrable system").
\item Sections \ref{QP-G-sss}, \ref{QP-Y-sss}:
a construction of  new explicit generators $(:Tr L^{[k]}(z):)$
\footnote{The basic definitions on  "normal ordering" $:...:$,
"critical level",
$U_{c=crit}(\widehat{gl_n})=U_{c=crit}(gl_n[t,t^{-1}],c)$, and
("Bethe") subalgebra in $U(gl_n[t]))$ are recalled in the Appendix
(section \ref{App-s})} in the center of
$U_{c=crit}(gl_n[t,t^{-1}],c)$ and $(Tr L^{[k]}(z))$ in the
commutative Bethe subalgebra in $U(gl_n[t])$; in the physics
language $(Tr L^{[k]}(z))$ are quantum conservation laws for the
quantum Gaudin system. Here $ L^{[k]}(z)$ are
 "quantum powers" defined in \CTKZ,
that "cures" some bad properties (i.e., failure of producing quantum
commuting operators) of the standard powers $L^{k}(z)$.
\footnote{E.g., $[Tr L_{gl_n-Gaudin}^{4}(z), Tr L_{gl_n~-Gaudin}^{2}(z)]\ne 0$
  \CRT\
this is expected to be a general phenomenon in the quantum case,
opposite to the classical one.} Similar results for the Yangian
case, and also  Newton and the Cayley-Hamilton identities for these
generators.
\item Section \ref{Nazar}: the construction of a further set of
explicit generators $(:Tr S^k L(z):)$ and  $(Tr S^k L(z))$ in the
center of $U_{c=crit}(gl_n[t,t^{-1}],c)$, and, respectively, in the
commutative Bethe subalgebra in $U(gl_n[t])$; in physical language,
$(Tr S^k L(z))$ are  quantum conservation laws for the quantum
Gaudin, Heisenberg, etc. integrable systems. We notice that some
results on the Cayley-Hamilton, Newton, MacMahon-Wronski relations
in $U(gl_n)$ by T. Umeda, M. Itoh \UmedaN,\UmedaW, \ItoN~
are implied by our constructions. These issues will be discussed
elsewhere.
\end{itemize}
Finally, we address the problem of  Sklyanin's separation of
variables problem; we present a conjectural construction of the
quantum separated variables for Yangian type systems, as well as
some remarks concerning the Gaudin case.

In the Appendix, we first collect definitions about the
center of $U_{c=crit}(gl_n[t,t^{-1}],c)$, that is, the commutative
Bethe subalgebra in $U(gl[t])$ (Section \ref{App-s}), and on the
relations of the present work with the geometric Langlands program
(Section \ref{GLCoC}).
In Sections \ref{Announce-s} and \ref{SV-s}
we give a new look and generalizations of the Capelli type identities {\Ok961,\MTVCap}, 
which appears to be related with the quantization of the Gaudin and Heisenberg's XXX systems
via the Wick ordering prescription.

Our starting point is the following remarkable construction
by D.Talalaev of the "quantum spectral curve" or the "quantum characteristic polynomial".
It solves the longstanding problems of the explicit efficient
construction of the center of $U_{c=crit}(gl_n[t,t^{-1}],c)$, the commutative Bethe subalgebra in $U(gl[t])$
\footnote{See \CTBig,
\cite{CM} } and explicitly produces a complete set of quantum
integrals of motion for the Gaudin system. It also has many other
applications (to the Bethe ansatz, separation of variables, the
Langlands correspondence, the Capelli identities, real algebraic
geometry etc). The right setup where the Talalaev's formula fits are
the ideas of E. K. Sklyanin on quantum separation of variables (see,
e.g., the surveys \SklA, \SklB).\\

{\em {\bf Theorem \T04}
Let  $L(z)$ be the Lax matrix of the $gl_n$-Gaudin model, and
consider the following differential operator in the variable $z$
(the Talalaev's  quantum characteristic polynomial or quantum
spectral curve):
\begin{equation}
\label{TT}
det^{column}(\partial_z - L_{gl_n-Gaudin}(z))= \sum_{i=0,...,n} QH_{i}(z) \partial_z^i,\\
\end{equation}
Then:
}
\[
\forall i,j\in 0,\ldots,n,\text{ and } u\in \CC,v\in \CC ~~~
[QH_{i}(z)|_{z=u} ,
QH_{j}(z)|_{z=v} ]=0.
\]
So taking $QH_i(z)$ for different $i,z=u\in \CC$ (or their residues
at poles, or other ``spectral invariant'' of the Lax
matrix $\partial_z - L_{gl_n-Gaudin}(z)$)
one obtains a full set of generators of quantum mutually
commuting conserved integrals of motion.\\
Actually, the theorem holds for all Lax matrices of ``$gl_n$-Gaudin
type (see definition \ref{GL-def}, section \ref{Gaud-Manin-ss}).\\
For the Yangian case D.Talalaev considers:
$det(1-e^{-\partial_z}T_{gl_n-Yangian}(z))$ (\T04 formula 9 page 6).
\paragraph{Remarks.}
1) In  \GLZ, \KonvalinkaA, \KonvalinkaB\
the names 
"right quantum matrices"; in
\RRT02, \RRTb05, \LT07~ the name ``left quantum group''
were used.
We prefer to use the name "\MMs".\\
2) Our case is different from the more general one
of  \GelfandGelfandRetakhWilson02,
where  {\em generic } matrices with noncommutative entries are
considered. In that case, there is no natural definition of the
determinant (rather, it was proposed to build linear algebra in
terms of elements $A^{-1}$, elements of $A^{-1}$ being called
"quasi-determinants", so that there are $n^2$ quasi-determinants
instead of a single notion of determinant in the commutative case).
Thus the analogues of linear algebra propositions are sometimes
quite different from the commutative case. Nevertheless some results
of the above mentioned authors can be fruitfully applied to
some questions here.\\
3) Those readers that are familiar with the $R$-matrix approach to
(quantum) integrable systems know that $L_{gl_n-Gaudin}(z)$ and
$T_{gl_n-Yangian}(z)$ satisfy quite {\em different} commutation
relations (\ref{Lax-G-com-rel},\ref{Yang-com-rel}) with "spectral
parameter", namely, in the first case we have linear $R$-matrix
commutation relations, while in the second one we have quadratic
relations. It is quite surprising that Talalaev's introduction of
$\p$ converts both to the {\em same} class of \MMs.

Moreover Manin's
relations do {\em not} contain  explicitly the {\em "spectral
parameter" $z$}.
Thus we can use simpler considerations (that is, without taking into
account the dependence on spectral parameter) and apply our results
also to the
theory $z$-dependent Lax/transfer matrices. Another feature of
insertion of $\p$ is that the somewhat {\em ad hoc shifts} in the
spectral parameter entering the formulas for the "$qdet(T(z))$"
known in the literature now {\em appear automatically} from the
usual column-determinant (see e.g. lemma \ref{ShiftsRem} page \pageref{ShiftsRem}).

All the considerations below work for an arbitrary field of
characteristic not equal to 2, but we restrict ourselves with $\CC$.

Subsequent paper(s)  \cite{CF} will contain more proofs, background material,
generalizations, e.g., to the (q)-super cases.

{\bf Acknowledgments} G.F. acknowledges support from the ESF
programme MISGAM, and the Marie Curie RTN ENIGMA (Contract number
MRTN-CT- 2004-5652). The work of A.C. has been partially supported
by the
Russian President Grant MK-5056.2007.1, the grant of Support for the
Scientific Schools 8004.2006.2, and by the INTAS grant
YSF-04-83-3396, by the ANR grant GIMP (Geometry and Integrability in
Mathematics and Physics). Part of work  was done during the visits
of AC to SISSA and the Milano-Bicocca University (under the INTAS
project), and to the University of Angers (under ANR grant GIMP).
A.C. is deeply indebted to SISSA and especially to B. Dubrovin, as
well as to the University of Angers and especially to V. Rubtsov,
for providing warm hospitality, excellent working conditions and
stimulating discussions. The authors are grateful to D. Talalaev, M.
Gekhtman, A. Molev,  E. Sklyanin, A. Silantiev  for multiple stimulating
discussions, to P. Pyatov for sharing with us his unpublished
results,
pointing out to the paper \GLZ ~
 and for multiple  stimulating discussions,
to S. Loktev for his friendly help, to A. Smirnov for Maple checks
of the conjectures and stimulating discussions, to Yu. Burman and N. Crampe
for their kind remarks.

\section{\MMs.
Definitions and elementary properties}
We herewith recall definitions and results of \ManinA, \cite{Manin,ManinBook91}
(with minor variations suited to the sequel of the paper). We first
remark that the notion (given in Definition \ref{D1}) of \MM\ with
elements in an associative ring ${\mathcal R}$ can be reformulated as
the condition that \[ \forall p,q,k,l~~[M_{pq}, M_{kl} ]= [M_{kq},
M_{pl} ]\>\text{ e.g. } [M_{11}, M_{22} ]= [M_{21}, M_{12} ],\>
[M_{11}, M_{2k} ]= [M_{21}, M_{1k} ]
\] Indeed, for
$q=l$  the above requirement yields that elements which
belong to the same column of $M$ commute among themselves. For
$q\neq l$, it is precisely the second condition of Definition~\ref{D1}.
{\Rem ~}
Precisely these relations were explicitly written in
\cite{Manin} (see chapter 6.1, especially Formula 1). Implicitly they
are contained in \ManinA -- the last sentence on page 198 contains a definition of the algebra
 $end(A)$ for an arbitrary quadratic Koszul algebra $A$. One can show (see the remarks
on the page 199 top) that $end(\CC[x_1,...,x_n])$ is the algebra generated by $M_{ij}$.

{\Rem ~\label{QGrem} Relation with quantum group $Fun_q(GL(n))$
(see \cite{Manin}, \cite{CF}).}
 Let us consider $n=2$. The q-analog of the relations
above is: $M_{21} M_{11}=q M_{11} M_{21}, M_{22} M_{12}=q M_{12} M_{22},
M_{11} M_{22} =  M_{22} M_{11} +q^{-1} M_{21} M_{12} -q M_{12} M_{21}$, this is precisely half
of the relations defining the quantum group $Fun_q(GL(2))$ (e.g. \ManinA~ page 192, formula 3).
This is true in general: q-\MMs~ are obtained by 
{\bf half} of the relations
of the corresponding quantum group $Fun_q(GL(n))$. Conversely
one can  define $Fun_q(GL(n))$ by matrices $M$ such that $M$ and $M^t$ are simultaneously q-\MMs~
\ManinA, \cite{Manin}.
In present paper we consider $q=1$ case.

The proposition below gives a more intrinsic characterization of
Manin's matrices 
(see \cite{Manin}, \cite{ManinBook91}:)

{\Prop\label{proposizione} Consider a rectangular $n\times m$-matrix
$M$, the polynomial algebra $\CC[x_1,...,x_m]$ and the Grassman
algebra $\CC[\psi_1,...,\psi_n]$ (i.e. $\psi_i^2=0, \psi_i \psi_j=
-\psi_j \psi_i$); suppose $x_i$ and $\psi_i$ commute with the matrix
elements $M_{p,q}$.
Consider the variables $\tilde x_i$, $\tilde
\psi_i$ defined by: \bea \left(\begin{array}{c}
\tilde x_1 \\
... \\
\tilde x_n
\end{array}\right)
=
\left(\begin{array}{ccc}
 M_{1,1} & ... &  M_{1,m} \\
... \\
 M_{n,1} & ... &  M_{n,m}
\end{array}\right)
\left(\begin{array}{c}
 x_1 \\
... \\
 x_m
\end{array}\right)
~~~~
(\tilde \psi_1, ... ,  \tilde \psi_m) =
(\psi_1, ... ,   \psi_n)
\left(\begin{array}{ccc}
 M_{1,1} & ... &  M_{1,m} \\
... \\
 M_{n,1} & ... &  M_{n,m}
\end{array}\right),
\eea that is the new variables are obtained via left action (in the
polynomial case) and right action (in the Grassmann case) of $M$ on
the old ones. Then the following three conditions are equivalent:
\begin{itemize}
\item The matrix $M$ is a \MM.
\item The variables
$\tilde x_i$ commute among themselves: $[\tilde x_i,\tilde x_j]=0$.
\item The variables $\tilde \psi_i$ anticommute among themselves:
$\tilde \psi_i\tilde \psi_j+\tilde \psi_j \tilde \psi_i =0$.
\end{itemize}}

{\bf Proof.} It is a straightforward calculation.

\vskip 0.4cm
Let us present some examples of \MMs.
{\Def \label{CF-def} A matrix $A$ with the elements in a noncommutative ring is
called a Cartier-Foata (see \cite{CF69,Fo79}) matrix if elements from different {\em rows} commute with
each other.}
{\Lem ~ Any Cartier-Foata matrix  is a 
\MM.
}\\
{\bf Proof}. The characteristic conditions of Definition \ref{D1} are trivially
satisfied in this case.\hfill\BX

The example below is related to our discussion of the
Capelli identities (section \ref{Cap-sssect})).
Let $x_{ij}, y_{ij}$ be commutative variables, and $X,Y$ be
$n\times k$ matrices with matrix elements $x_{ij}, y_{ij}$. Also, let
$\partial_X, \partial_Y$ be  $n\times k$ matrices with matrix
elements $\frac{\partial}{ \partial x_{i,j}}, \frac{\partial}{
\partial y_{i,j}}$. Let $z$ be a variable commuting with $y_{ij}$.
The following $2n\times 2k$, and $(n+k)\times (n+k)$ matrices, with
elements in the ring of differential operators in the variables
$X,Y$ are easily seen to be
\MMs: \begin{equation}\label{MTVex}
 \left(\begin{array}{cc}
X  & \partial_Y \\
Y & \partial_X \\
\end{array}\right),\quad\text{and}\quad
\left(\begin{array}{cc}
z ~ 1_{k\times k}  & (\partial_Y)^t \\
Y & \partial_z ~ 1_{n\times n}\\
\end{array}\right).
\end{equation}

\subsection{The determinant} \label{Det-ss}
{\Def
Let $M$ be a \MM.
Define the determinant of $M$ by  column expansion: \bea det
M=det^{column} M=\sum_{\sigma\in S_n} (-1)^\sigma
\prod^{\curvearrowright}_{i=1,...,n} M_{\sigma(i),i}, \eea where
$S_n$ is the group of
permutations of $n$ letters, and the
symbol $\curvearrowright$ means that in the product
$\prod_{i=1,...,n} M_{\sigma(i),i}$ one writes at first the elements
from the first column, then from the second column and so on and so
forth.}

{\Lem \label{det-Col-indep} The determinant of a \MM~ does not depend
on the order of the columns in the column expansion, i.e.,
\bea \forall p\in S_n ~~~ det^{column} M=\sum_{\sigma\in S_n}
(-1)^\sigma \prod^{\curvearrowright}_{i=1,...,n}
M_{\sigma(p(i)),p(i)} \eea }\noindent  {\bf Proof}. Since any
permutation can be
presented as a product of transpositions of neighbors $(i,i+1)$ it
is enough to prove the proposition for such transpositions. But for
them it follows from the equality of the commutators of cross
elements (formula \ref{Cross-term-propert}). \hfill\BX

{\Ex ~} For the case $n$=2, we have
\bea
det^{col}\left(\begin{array}{cc}
a & b\\
c & d \\
\end{array}\right)
\stackrel{def}{=} ad-cb\stackrel{lemma}{=}da-bc.
\eea Indeed, this is a restatement of the second condition of
Definition \ref{D1}.
\subsubsection{Elementary properties} The
following properties are simple consequences of the definition of
\MM.

\begin{enumerate}
{\item ~ Any matrix with commuting elements is a \MM. } \vskip 0mm
\item ~ Any submatrix of a \GT\  matrix is again a \GT\ matrix.
\item ~ If $A$, $B$ are  \GT\  matrices and $\forall i,j,k,l:
[A_{ij},B_{kl}]=0$, then $A+B$  is again a \GT\ matrix.
{\item ~ If $A$ is a \GT\  matrix, $c$ is constant,  then $cA$ is a \MM. } \vskip 0mm
{\item ~ If $A$ is a \GT\  matrix, $C$ is constant matrix,  then $CA$
and $AC$ are \GT\ matrices and $det(CA)=det(AC)=det(C)det(A)$. }
\vskip 0mm
{\item ~ If $A,B$ are  \GT\  matrices and $\forall i,j,k,l:
[A_{ij},B_{kl}]=0$, then
 $AB$ is a \GT\ matrix and $det(AB)=det(A)det(B)$. }
\vskip 0mm
{\item ~ If $A$ is a
\GT\  matrix, then one can exchange the $i$-th and the $j$-th
columns(rows);
 one can put
  $i$-th column(row) on $j$-th place (erasing $j$-th column(row));
one can add new column(row) to matrix $A$ which is equal to one of
the columns(rows) of the matrix $A$; one can add the $i$-th
column(row) multiplied by any constant to the $j$-th  column(row);
in all cases the resulting matrix will be again a
\GT\   matrix.
}
\vskip 0mm
{\item ~ If $A$ and simultaneously $A^t$  are
\GT\    matrices,
then all elements $A_{i,j}$ commute with each other. (A q-analog of this lemma
says that if  $A$ and simultaneously $A^t$  are
q-\GT\ , then $A$ is quantum matrix: $"R\one A \two A=\two A \one A R"$
(\ManinA, \cite{Manin} 
)).}
\item
The exchange of two columns in a \MM~ changes the  sign of the determinant. If
 two columns or two rows in a 
\MM~ $M$ coincide,  then $det(M)=0$.\\
\PRF The first assertion is obvious. For the second, (in the case of
equal $i$-th and $j$-th column) use the column expansion of $det(M)$,
taking any permutation with first two elements $i$ and $j$. One sees
that $det(M)$ is a sum of elements of the form $(xy-yx)(z)$, where
$x,y$ belong to the same column. Using the column commutativity
$[x,y]=0$ we get the result. For the case of two coinciding rows the
assertion follows at once, without the help of column
commutativity. \hfill\BX {\item Since any submatrix of a \GT ~ matrix is a
\GT\ matrix one has a natural definition of minors and again one can
choose an arbitrary order of columns (rows) to define their
determinant.}
\end{enumerate}

\subsubsection{The permanent} The permanent of a matrix is a
multilinear function of its entries,
similar to the determinant, with no $(-1)^{sign (\sigma)}$ in its
definition. We will make use of it in section \ref{MW-s} below.
{\Def
Let $M$ be a \MM.
Let us define the permanent of $M$ by  {\bf row} expansion
\footnote{Remark the difference with the definition of the determinant,
where one uses column expansion} 
: \bea \label{M-perm-f} perm M=perm^{row} M=\sum_{\sigma\in S_n}
\prod^\curvearrowright_{i=1,...,n} M_{i,\sigma(i)}, \eea
}
{\Lem \label{perm-Col-indep} The permanent of \MMs~
does not depend on the order of rows in the row expansion:
\vskip -10mm \bea \forall p\in S_n ~~~ perm^{row} M=\sum_{\sigma\in
S_n}  \prod^\curvearrowright_{i=1,...,n} M_{p(i),\sigma(p(i))} \eea
} {\bf Proof}. Since any permutation can be presented as a product of
transpositions of neighbors $(i,i+1)$ it is enough to prove the
proposition for such transpositions. But for them it follows from
the equality of the commutators of cross elements (formula
\ref{Cross-term-propert}).
\hfill\BX
{\Ex ~}
\vskip -10mm
\bea
perm^{row}\left(\begin{array}{cc}
a & b\\
c & d \\
\end{array}\right)
\stackrel{def}{=} ad+bc \stackrel{lemma}{=}da+cb,
\eea

{\Lem \label{perm-Row-indep} It is easy to see that the row-permanent of
an arbitrary matrix $M$ (even without conditions of commutativity)
does not change  under any permutation of  columns.
}

\section{Lax matrices as  \MMs }\label{sect3}
In the modern theory of integrable systems, "Lax matrices" are also
called "Lax operators", "transfer matrices", "monodromy matrices"
for historical reasons. It is outside of the aim of the paper to make
an attempt to provide a kind of formal definition of a Lax matrix
for an integrable system. However, it is fair to say that there is a
set of properties expected from a "good" Lax matrix;  let us now
recall those  property relevant for our exposition.

Let $M$ be a symplectic manifold and $H_i$ be a set of
Poisson-commuting functions defining an integrable system on it (for
simplicity, let us consider the Hamiltonian flow $X_1$ of one of
these functions, say $H_1$)
A {\em Lax matrix $L(z)$} for 
$(M,H_i)$ is a matrix, whose matrix elements are functions on $M$,
possibly depending \footnote{let us assume for simplicity of
presentation
 that $ L(z)$ is just a polynomial function of the formal  parameter $z$}
on a formal parameter $z$, satisfying the following two
characteristic properties:
\begin{enumerate}
\item
The evolution of $L(z)$ along $X_1$ is of the form
\[
\dot L(z)=[M_1(z),L(z)].
\]
\item
{\em The characteristic polynomial} 
$det(\l - L(z))=\sum_{i,j} H_{i,j} z^j\l^i$
produces "all Liouville" integrals of motion, i.e.
\begin{itemize}
\item $\forall i, j ~~ H_{i,j}$ Poisson commute among themselves and with the
given functions $H_i$ \bea \forall k, i, j ~~ \{ H_k , H_{i,j} \}
=0, ~~~~~~~~~~~~~~ \forall i,j, k,l ~~ \{ H_{i,j} , H_{k,l} \} =0
\eea
\item
 All $H_i$ can be expressed via $H_{k,l}$ and vice versa.
\end{itemize}
\end{enumerate}
Poisson algebras of commutative functions on  manifolds $M$ are
related to classical mechanics. As it is well known, in quantum
mechanics one is lead to consider a family of
 non commutative,
but associative algebras  $\widehat{Fun(M)}_{\hbar}$; $\hbar$ is a
formal parameter in mathematics and Plank's constant in physics. For
$\hbar=0$ the algebra $\widehat{Fun(M)}_{\hbar}$ coincide with the
commutative  algebra $Fun(M)$ of functions on a manifold $M$, as
well as  Poisson brackets are related to commutators: $[f,g]=\hbar
\{f,g \} mod (\hbar^2)$, for $f,g\in \widehat{Fun(M)}_{\hbar}$.
(See, e.g.,  \Kontsevich97).

The standard example is the algebra of functions on $\CC^{2n}$ -
$\CC[p_i,q_i]$ with the Poisson bracket $\{p_i, q_j\}=\delta_{i,j},
\{p_i, p_j\}= \{q_i, q_j\}=0$  whose quantization is the Heisenberg
(Weyl) algebra generated by $\hat p_i, \hat q_i$ and the relations
$[\hat p_i, \hat q_j]=\hbar \delta_{i,j}, [\hat p_i,\hat p_j]=[\hat
q_i,\hat q_j]=0$. (We will usually put $\hbar=1$).

Within this framework, on can look  for a {\em Quantum Lax matrix}
for a quantum systems as a matrix satisfying
\begin{itemize}
\item $\hat L(z)$
is a matrix which matrix elements are elements from
$\widehat{Fun(M)}_\hbar$, usually depending on a formal parameter $z$
(thus, a  quantum Lax matrix is a matrix with {\em noncommutative }
elements).
\item in the  classical limit $\hbar \to 0$ one has  $ \hat L(z) \to L(z)^{classical}$.
\end{itemize}
It is quite natural to look, in the integrable case, for  a kind of
determinantal formula: $"det(\hat \l - \hat L(\hat z))" =\sum_k \hat
\l^k \hat H_k (\hat z)$ to produce quantum integrals of motion:
$[\hat H_k(z), \hat H_l (u)]=0$,
and, possibly,
to satisfy other important properties (see \CTBig).

\subsection{The Gaudin case \label{Gaud-Manin-ss} }
 The algebra of symmetries of the $gl_n$-Gaudin
(\cite{Gaudin76,GaudinBook}) integrable system is the Lie algebra
$gl[t]$. The Lax matrix for the Gaudin system
 is a convenient way to combine  generators of $gl[t]$ (or its factor algebras)
 into one generating matrix-valued function
(see formula \ref{Lax-glt}b~ below for $gl[t]$ itself). Quantum
commuting Hamiltonians (quantum conservation laws) arise from the
maximal commutative Bethe subalgebra in $U(gl_n[t])$. Explicit
efficient expressions for generators were first obtained by D.
Talalaev. This subalgebra is actually an image of the center of
$U_{c=crit}(gl_n[t,t^{-1}],c)$ under the natural projection
$U_{c=crit}(gl_n[t,t^{-1}],c) \to U(gl_n[t])$ (see Appendix, section
\ref{App-s}). This explains the mathematical meaning and the
importance of the subject.

{\Def \label{GL-def} Let $R$ be an associative algebra over $\CC$.
Let us call a matrix $L(z)$ with elements in $R((z))$ (i.e. $L(z)\in Mat_n\otimes R\otimes \CC ((z))$)
\footnote{ or $L(z)\in Mat_n\otimes R\otimes \CC ((z^{-1}))$; in many cases $L(z)$ is just a rational function of $z$}
a {\bf  Lax matrix of $gl_n$-Gaudin type} iff:}
\bea \label{Lax-G-com-rel}
&& [ L_{ij}(u), L_{kl}(v) ]=\frac{1}{u-v}(L_{il}(v)\delta_{jk}-
L_{il}(u)\delta_{jk} - L_{kj}(v)\delta_{li}+  L_{kj}(u)\delta_{li}).
\eea
More precisely this is "rational" $gl_n$-Gaudin type; as it is
well-known,there are trigonometric and elliptic versions as well as
generalizations to semisimple Lie (super)-algebras.
\\
We recall  the matrix ("Leningrad's") notations:  let $\Pi\in
Mat_n\otimes Mat_n$ be the  permutation matrix: $\Pi (a\otimes b) =
b \otimes a$, and consider $L(z)\otimes 1, 1 \otimes L(u) \in
Mat_n\otimes Mat_n\otimes R \otimes \CC((z))\otimes \CC((u))$.
Formula \ref{Lax-G-com-rel} can be compactly written as follows:
\bea \label{Lax-G-com-rel-matrix}
&& [ L(z)\otimes 1, 1 \otimes L(u)] = [\frac{\Pi}{z-u}, L(z)\otimes 1 + 1 \otimes L(u)],
\eea
that is,
a Lax matrix of $gl_n$-Gaudin type is a Lax matrix with {\em linear}
$r$-matrix structure (the $r$-matrix being $r=\frac{\Pi}{z-u}$).
{\Rem From(\ref{Lax-G-com-rel}) it follows that
\[
[ L_{ij}(u), L_{kl}(u) ]=-(\partial_u L_{il}(u)\delta_{jk} -
\partial_u L_{kj}(u)\delta_{li}).
\]}
Here is our {\em first main observation}:
{\Prop \label{G-M-pr}
Consider $\p \pm L(z)\in Mat_n\otimes R \otimes \CC((z))[\p]$,
where $L(z)$ is a Lax
matrix of the $gl_n$-Gaudin type
above; then:
\bea \label{G-M-fml}
(\p-L(z)), ~~~ (\p+L(z))^t \mbox{ ~~  are \MMs}
\eea } \PRF The proof is a straightforward computation.

{\Ex ~} Let us show this property in the example of $2\times 2$
case:
\begin{equation}
\begin{split}
&\p-L(z) = \left(\begin{array}{cc}
\p-L_{11}(z) &  -L_{12}(z) \\
-L_{21}(z) & \p-L_{22}(z) \\
\end{array}\right), \\
&\text{Column $1$  commutativity}:~~ [\p-L_{11}(z), -L_{21}(z)]=-L_{21}'(z)+ [L_{11}(z), L_{21}(z)]=0,\\
&\text{Cross term relation} : ~~[\p-L_{11}(z),
\p-L_{22}(z)]=-L_{22}'(z)+ L_{11}'(z)= [L_{21}(z), L_{12}(z)].
\end{split}
\end{equation}

The following well-known fact easily follows from \ref{Lax-G-com-rel-matrix}:
{\Prop Let $L(z)$ and $\tilde L(z)$ be Lax matrices of the Gaudin type with pairwise commuting elements
(i.e. $\forall i,j,k,l,z,u:~~[L_{ij}(z), \tilde L_{kl}(u)]=0$)
then $L(z) + \tilde L(z)$ is again a Lax matrices of the  Gaudin type.}

{\Rem ~} The classical counterpart of the commutation relations
(\ref{Lax-G-com-rel}) -- namely, with Poisson brackets replacing
commutators) and of the corresponding Lax matrices are associated with
a huge number of integrable systems. In particular, among integrable
systems leading to Lax representations of Gaudin type we can
mention: the Neumann's geodesic motion on ellipsoids and the
Neumann-Rosochatius systems, the $gl_n$-Manakov-Mishenko-Fomenko tops
(see, e.g., \AHH96)\
to the $so(3)$ case);
the Lagrange top and the Clebsh system (see, e.g., \PetreraRagnisco07);
"polynomial matrices" (and so finite gap solutions of
KdV-type equations, although non-trivial
reductions of the $r$-matrix bracket have to be
considered) \cite{Mumford84}, \ReymanSemenov79,
\AdlerMoerbeke80, \Beaville90;
the general systems studied in \AdamsHarnadPreviato88;
the Jaynes-Cummings-Gaudin system, \BT;~
as well as
Hitchin's system on singular rational curves \CTHB.~
Bending flows, Gelfand-Ceitlin, Jucys-Murphy integrable
systems are degenerations of the Gaudin system (\FMA,\FMB,\CFR).

Moreover, Schlesinger system of isomonodromy deformations and
Knizhnik-Zamolodchikov equation are
 non-autonomous versions of the classical and, respectively,
 quantum Gaudin Systems \DubrovinM2,  
\Reshetikhin92. 
It is well-known that Nahm's monopole equations can be rewritten as
a Gaudin system (see e.g. \GGM97 
(section 6.1 page 24), where applications to Seiberg-Witten theory
of N=2 supersymmetric Yang-Mills are also discussed).
In the influential paper \FFR ~
(see, also the survey \Frenkel95) 
the relation with the geometric Langlands correspondence was
discovered. It might also be noticed that a great number of systems
in condensed matter theory are related to Gaudin-type systems, see,
e.g., \OrtizSommaDukelskyRombouts04.

In this paper, we shall concentrate on the following
few examples of Gaudin type Lax matrices.

\paragraph{The simplest example}
Let $K$ be an arbitrary constant matrix, and $n,k\in {\mathbb N}$,
and $z_1,\ldots,z_k$ arbitrary points in the complex plane. Consider
\bea \label{Lax-G-ex1}
L(z) = K+ \sum_{i=1,...,\mpn} \frac{1}{z-z_i}
\left(\begin{array}{c}
\hat q_{1,i} \\
... \\
\hat q_{\ms,i } \\
\end{array}\right)
\left(\begin{array}{cccc} \hat p_{1,i} &  ... &  \hat p_{\ms,i}
\end{array}\right) = K+ \hat Q~
diag(\frac{1}{(z-z_1)},...,\frac{1}{(z-z_\mpn)} )~ \hat P^t \eea
where $\hat p_{i,j},\hat q_{i,j}$,  $i=1,...,\ms;j=1,...,\mpn$ are
the standard generators of the standard Heisenberg algebra $[ \hat
p_{i,j}, \hat q_{k,l} ]=\delta_{i,k} \delta_{j,l}$, $[ \hat p_{i,j},
\hat  p_{k,l} ]= [ \hat q_{i,j},\hat  q_{k,l} ]=0$. Also, $\hat
Q,\hat P$ are $n\times k$-rectangular matrices with elements $\hat
Q_{i,j}=\hat q_{i,k}$, $\hat P_{i,j}=\hat p_{i,j}$.

One can see that this Lax matrix satisfies relations  \ref{Lax-G-com-rel-matrix}
(\ref{Lax-G-com-rel}),
and so by proposition \ref{G-M-pr} above $(\p - L(z)), (\p + L(z))^t $ are \MMs.


\paragraph{The standard example}
Consider $gl_n\oplus ... \oplus gl_n$ and denote by $e_{kl}^i$ the
standard basis element from the $i$-th copy of the direct sum
$gl_n\oplus ... \oplus gl_n$.
The standard Lax matrix for the Gaudin system is:
\bea \label{Lax-st}
L_{gl_n-Gaudin~standard}(z) =
\sum_{i=1,...,\mpn} \frac{1}{z-z_i}
\left(\begin{array}{ccc}
 e_{1,1}^{i} &  ... & e_{1,n}^{i}  \\
 ... & ... &...  \\
 e_{n,1}^{i}  & ... & e_{n,n}^{i}
\end{array}\right)
\eea
$L_{gl_n-Gaudin~standard}(z)\in Mat_n \otimes U(gl_n\oplus ... \oplus gl_n)\otimes \CC(z).$

One can see that this Lax matrix satisfies relations  \ref{Lax-G-com-rel-matrix} (\ref{Lax-G-com-rel}),
and so by proposition \ref{G-M-pr} above $\p - L_{gl_n-Gaudin~standard} (z)$,
$(\p + L_{gl_n-Gaudin~standard} (z))^t$ are \MMs.

\paragraph{The $gl_n$ and $gl_n[t]$ examples.}
Consider the Lie algebra $gl_n$, with $e_{ij}$ its standard linear
basis;  consider the polynomial Lie algebra $gl_n[t]/t^\mpn$,
$\mpn=1,...,\infty$. The Lax matrices for $gl_n$ and
$gl_n[t]/t^\mpn$ are the following: \bea \label{Lax-gl} L_{gl_n}(z)
= \frac{1}{z} \left(\begin{array}{ccc}
 e_{1,1} &  ... & e_{1,n}  \\
 ... & ... &...  \\
 e_{n,1}  & ... & e_{n,n}
\end{array}\right)
~~~~~~ \label{Lax-glt}
L_{gl_n[t]}(z) = \sum_{i=1,...,\mpn} \frac{1}{z^i}
\left(\begin{array}{ccc}
 e_{1,1}t^{i-1} &  ... & e_{1,n}t^{i-1}  \\
 ... & ... &...  \\
 e_{n,1}t^{i-1}  & ... & e_{n,n}t^{i-1}
\end{array}\right).
\eea
$L_{gl_n}(z)\in Mat_n \otimes U(gl_n)\otimes \CC(z)$; ~~
$L_{gl_n[t]}(z)\in Mat_n \otimes U(gl_n[t])\otimes \CC((z^{-1}))$, where $U(g)$ is the universal enveloping
algebra of $g$. Note: $zL_{gl_n}(z)$ coincide with the expression (12) section 2.1 \Kir.

One can see that these Lax matrices satisfy relations
\ref{Lax-G-com-rel-matrix} (\ref{Lax-G-com-rel}), and so by
proposition \ref{G-M-pr} above $\p - L_{gl_n} (z), \p -L_{gl_n[t]}(z) $,
$(\p + L_{gl_n} (z))^t, (\p +L_{gl_n[t]} (z))^t $  are \MMs.

\subsection{The Yangian (Heisenberg chain) case.}
\label{Yang-ss} The Algebra of symmetries of the Heisenberg XXX spin
chain, the quantum  Toda system and several other integrable systems
is a Hopf algebra called Yangian. It was implicitly defined in the
works of Faddeev's school; a concise mathematical treatment and deep
results were given in \cite{Dr85}
(see also the excellent more recent surveys \Molev02,\MNO).
It is a deformation of the universal enveloping algebra of $gl[t]$.
A Lax matrix of Yangian type  is a convenient way to combine
generators of the Yangian (or its factor algebras) into one
generating matrix-valued function. Following standard notations we
will write $T(z)$ ("Transfer matrix") instead of $L(z)$ in the case
of the Yangian type Lax matrices.

{\Def \label{Yang-def} Let $\mathcal R$ be an associative algebra
over $\CC$. Let us call a matrix $T(z)$ with elements in ${\mathcal
R}((z))$ (i.e. $T(z)\in Mat_n\otimes {\mathcal R}\otimes \CC ((z))$)
\footnote{ or $T(z)\in Mat_n\otimes {\mathcal R}\otimes \CC
((z^{-1}))$; in many cases $T(z)$ is just a rational function of
$z$}
a {\bf  Lax matrix of Yangian type} iff:}
\bea \label{Yang-com-rel}
&& [ T_{ij}(u), T_{kl}(v) ]=\frac{1}{u-v}(T_{kj}(u)T_{il}(v)-T_{kj}(v)T_{il}(u))
\eea
(See \Molev02 ~
page 4 formula 2.3.
More precisely this is the case of the $gl_n$-Yangian, there are
generalizations to the semisimple Lie (twisted)-(super)-algebras).

In matrix (Leningrad's) notations we have, with $\Pi\in Mat_n\otimes
Mat_n$ the permutation matrix: $\Pi (a\otimes b) =  b \otimes a$,
and $T(z)\otimes 1, 1 \otimes T(u) \in  Mat_n\otimes Mat_n\otimes
{\mathcal R} \otimes \CC((z))\otimes \CC((u))$, $R(z-u)=(1\otimes
1-\frac{\Pi}{z-u})$ that formula \ref{Lax-G-com-rel} can be written
as follows, as a {\em quadratic} $R$-matrix relation:
(see \Molev02 ~
page 6 proposition 2.3 formula 2.14).
\bea \label{Yang-com-rel-matrix}
  (1\otimes 1-\frac{\Pi}{z-u})  (T(z)\otimes 1) ~ (1 \otimes T(u))  =
 (1 \otimes T(u)) ~ (T(z)\otimes 1)  (1\otimes 1-\frac{\Pi}{z-u}) \\
\mbox{ Or shortly: } R(z-u) \one T(z) \two T(u)= \two T(u) \one T(z) R(z-u).
\eea
Our {\em second main observation} is:
{\Prop ~ \label{Yang-GT} If
$T(z)$ is a Lax matrix of the Yangian type then
\bea T(z)e^{-\p},  ~~ (e^{\p} T(z))^t  ~~~~~~ \mbox{ are  \MMs.
}
\eea
Here $e^{\pm \p} T(z) \in Mat_n \otimes Y(gl_n)\otimes
\CC[[1/z,e^{-\p}]]$,
 (see definition
\ref{Yang-def} above). Actually $\forall c\in \CC: T(z+c)$ is also of the Yangian type, so
$T(z+c)e^{-\p}= e^{-\p} T(z+c+1),  (e^{\p} T(z+c))^t=( T(z+c+1)e^{\p})^t$ are \MMs.
\footnote{Observe an analogy. It is known (\cite{Manin}) that,
if $M$ is matrix defining $Fun_q(GL_n)$ i.e. it satisfies
$R\one M \two M=\two M \one M R$ {\em without } spectral parameter,
then $M$ and simultaneously $M^t$ are q-\MMs.
Here we have: $T(z)$ satisfies $R(z-u) \one T(z) \two T(u) =\two T(u) \one T(z) R(z-u)$
{\em with } spectral parameter and proposition claims that $e^{-\p} T(z)  $
and $ (e^{\p} T(z))^t  $  are simultaneously \MMs. Both cases are related to $GL_n$,
extension to the other Lie groups is not known for the moment.
}
}
\\
\PRF As in the Gaudin-type case, it follows from a straightforward
computation. Actually it is known that any $k\times k$ submatrix of the Yangian type matrix
is again a Yangian type matrix of $k\times k$-size, so essentially it is necessary to check only $2\times 2$ case.
\\
The following well-known fact easily follows from
\ref{Yang-com-rel-matrix}:
{\Prop Let $T(z)$ and $\tilde T(z)$ be
Lax matrices of Yangian type with pairwise commuting elements
(i.e. $\forall i,j,k,l,z,u:~~[T_{ij}(z), \tilde T_{kl}(u)]=0$); then
the product $T(z)\tilde T(z)$ is again a Lax matrix of Yangian
type.}
{\Lem ~\label{ShiftsRem} } Observe that $det^{column} ( T(z)e^{-\p})= qdet(T(z)) e^{-n\p}$,
where\\ $qdet(M(z))=\sum_{\sigma\in S_n} \prod_k
M_{\sigma(k),k}(z-k+1)$ (see e.g. formula 2.24 page 10 \Molev02).
From 2.25 loc.cit. one sees:
$det^{row} (e^{\p} T(z))= e^{n\p} qdet(T(z)) $.
So "shifts" in $z$ in formulas for $qdet$ now appear automatically from the usual column/row determinants.
\\
\vskip 0.5cm
Let us herewith list a couple of remarkable examples of Yangian-type
Lax matrices.
\paragraph{The Toda system}
Consider the Heisenberg algebra generated by $\hat p_i, \hat q_i$,
$i=1,...,n$ and relations $[ \hat p_i, \hat q_j ]=\delta_{i,j}$,
$[\hat p_i, \hat p_j ]= [ \hat q_i, \hat q_j ]=0$.
Define
\bea\label{ln} T_{Toda}(z)\,=\, \prod_{i=1,...,n}
\left(\begin{array}{cc}
z-\hat p_i & e^{-\hat q_i}\\
-e^{\hat q_i} & 0
\end{array}\right)
\eea
One can see that this Lax matrix satisfies relations  \ref{Yang-com-rel-matrix} (\ref{Yang-com-rel}),
and so by proposition \ref{Yang-GT} above $e^{-\p} T_{Toda}(z)$ is a
\MM. One can easily see that, with the identification $q_{n+1}=q_1$,
the coefficient $C_{n-1}$ of $z^{n-2}$ of $Tr T_{Toda}(z)$ equals
$\sum_{k<l\le n} \hat p_k \hat p_l - \sum_{i=1,...,n} e^{\hat
q_i-\hat q_{i+1}}$, and the coefficient $C_{n-1}$ of $z^{n-1}$ is
$-\sum_{i=1}^n p_i$. Thus $\frac12 C_{n-1}^2-C_{n-2}$ is the
physical Hamiltonian (i.e., the energy) of the periodic Toda chain.

\paragraph{The Heisenberg's XXX system}
Consider the Heisenberg algebra generated by
$[ \hat p_{i,j}, \hat q_{k,l} ]=\delta_{i,k} \delta_{j,l}$,
$[ \hat  p_{i,j}, \hat  p_{k,l} ]= [ \hat q_{i,j},\hat  q_{k,l} ]=0$;
the quantum Lax matrix for the simplest case of the $gl_n$-Heisenberg spin chain integrable systems
can be given as follows:
\bea
\label{Lax-XXX-1}
T_{XXX~"simplest"}(z) = \prod_{i=1,...,k} \left( 1_{n\times n}
+\frac{1}{z-z_i} \left(\begin{array}{c}
\hat q_{1,i} \\
... \\
\hat q_{n,i } \\
\end{array}\right)
\left(\begin{array}{cccc}
\hat p_{1,i} &  ... &  \hat p_{n,i} \end{array}\right)
\right)
\eea
One can see that this Lax matrix satisfies relations  \ref{Yang-com-rel-matrix} (\ref{Yang-com-rel}),
and so by proposition \ref{Yang-GT} above $e^{-\p} T_{XXX}(z)$ is a
\MM.

Consider $gl_n\oplus ... \oplus gl_n$ and denote by $e_{kl}^i$ the
standard basis element of the $i$-th copy of the direct sum
$gl_n\oplus ... \oplus gl_n$.
The standard Lax matrix for the $GL(N)$-Heisenberg's XXX system is:
\bea \label{Lax-st-XXX}
T_{gl_n-XXX-standard}(z) =
\prod_{i=1,...,\mpn} \left( 1_{n\times n}+  \frac{1}{z-z_i}
\left(\begin{array}{ccc}
 e_{1,1}^{i} &  ... & e_{1,n}^{i}  \\
 ... & ... &...  \\
 e_{n,1}^{i}  & ... & e_{n,n}^{i}
\end{array}\right) \right)
\eea
$T_{gl_n-XXX-standard}(z)\in Mat_n \otimes U(gl_n\oplus ... \oplus gl_n)\otimes \CC(z)$
.

One can see that this Lax matrix satisfies the relations
\ref{Yang-com-rel-matrix} (\ref{Yang-com-rel}),
and so by proposition \ref{Yang-GT} above $e^{-\p} T_{gl_n-XXX-standard}(z)$
is a \MM.

\section{ Algebraic properties of \MMs\  and their applications.}
In this section we will derive a few less elementary
properties of \MM, and give
applications thereof to integrable systems.
\subsection{ Cramer's formula \label{Cramer-s}}
{\Prop\label{ad-inv} \cite{Manin} Let $M$ be a
\GT\  matrix and denote by $M^{adj}$ the adjoint matrix defined in
the standard way, (i.e.
 $M^{adj}_{kl}=(-1)^{k+l}det^{column}(\widehat{M}_{l k})$ where
$\widehat{M}_{l k}$ is the $(n-1)\times (n-1)$ submatrix of $M$ obtained
removing the l-th row and the k-th column.
Then the same formula as in the commutative case holds true, that is,
\begin{equation}\label{Ead}
 M^{adj} M= det^{column}(M) ~Id
\end{equation}
If $M^t$ is a \MM, then  $M^{adj}$ is defined by row-determinants
and $M M^{adj} = det^{column}(M^t) ~Id=det^{row}(M) ~Id$.} {\Ex ~}
In the $2\times 2$ case we have: \bea && \left(\begin{array}{cc}
d & -b \\
-c & a
\end{array}\right)
\left(\begin{array}{cc}
a & b \\
c & d
\end{array}\right)
=
\left(\begin{array}{cc}
da-bc & db-bd \\
-ca+ac & -cb+ad
\end{array}\right)
=
\left(\begin{array}{cc}
ad-cb & 0 \\
0 & ad-cb
\end{array}\right),
\eea where the characteristic commutation relations of a \MM\ have
been taken into account.

\PRF One can see that the equality: $\forall~i:(M^{adj}
M)_{i,i}=det^{col}(M)$, follows from the fact that $det^{col}(M)$
does not depend on the order of the column expansion of the
determinant. This independence was proved above (lemma
\ref{det-Col-indep}).
Let us introduce a matrix $\tilde M$ as follows. Take the matrix $M$
and set the $i$-th column equals to the $j$-th column; denote the
resulting matrix by $\tilde M$.
Note that $det^{col}(\tilde M)=0$ precisely gives $(M^{adj}
M)_{i,j}=0$ for $i\ne j$. To prove that $det^{col}(\tilde M)=0$ we
argue as follows. Clearly $\tilde M$ is a \MM. Lemma
\ref{det-Col-indep} allows to calculate the determinant taking the
elements first from $i$-th column, then $j$-th, then other elements
from the other columns. Now it is quite clear that $det^{col}(\tilde
M)=0$, due to it is the sum of the elements of the form
$(xy-yx)(z)=0$, where $x,y$ are the elements from the $i$-th and
$j$-th of $\tilde M$, so from $j$-th column of $M$. By Manin's
property elements from the same column of a matrix $M$ commute, so
$xy-yx=0$, so $det^{col}(\tilde M)=0$.
\hfill\BX
{\Rem ~} The only difference with the commutative case is that, in the
equality (\ref{Ead}) the order of the products of $M^{adj}$ and $M$ has to be
kept in mind.

\subsubsection{Application to the Knizhnik-Zamolodchikov equation
\label{KZ}}

We can give a very simple proof of the formula relating the
solutions of  KZ with coupling constant $\kappa=1$ to the solutions
of the equation defined by Talalaev's formula:
\[
det(\p -L(z))
Q(z)=0.
\]
This result was first obtained in \CTGoper~ 
in a more complicated manner (see also \CTKZ).~

The standard  $gl_n$-KZ-equation (actually a system of equations)  \KZ~ 
for an unknown function  $\Psi(z_0, z_1,..., z_\mpn)$ is the following.
Let $V_{i}, i=0,...,\mpn$ be representations of $gl_n$,   $\Psi(z_0, z_1,..., z_\mpn)$
is\\ $V_{0}\otimes V_1 \otimes ... \otimes V_\mpn$ valued function of $z_p$, $p=0,...,\mpn$.
Let $e_{a b}$ be the standard basis in $gl_n$, $e_{a b}^i$  the corresponding basis in
$gl_n\oplus ... \oplus gl_n$ (upper index $i$ is a number of copy of $gl_n$).
Let $\pi$ denote the representation of $gl_n\oplus ... \oplus gl_n$ in $V_{0}\otimes V_1 \otimes ... \otimes V_\mpn$.
\bea
\mbox{KZ system. ~~~~~~~~~ j-th equation: } ~~~
\left(\partial_{z_j} -  \kappa \sum_{i=1...\mpn}
\frac{  \sum_{ab}\pi(e_{ab}^j) \otimes \pi(e_{ab}^{(i)})}
{z_j-z_i} \right)\Psi(z) =0.
\eea
$j=0,...,\mpn$, and $\kappa\in \CC$ is some constant.

Let us restrict to the case $V_0$ is $\CC^n$ and denote $z_0$ by $z$.
Let us now consider only one equation, say the first and observe the following relation with the Gaudin Lax matrix:
\bea \label{KZ-fml}
\left(\p-  \kappa \sum_{i=1...\mpn}
\frac{  \sum_{ab}E_{ab}\otimes \pi_i(e_{ab}^{(i)})}
{z-z_i} \right)\Psi(z)=\pi(\p - \kappa L_{Gaudin}(z))\Psi(z)=0, ~~
\Psi(z)
=\left(
\begin{array}{c}
\Psi_1(z) \\
... \\
\Psi_n(z)
\end{array}\right)
\eea
So we came to:\\
{\bf Observation.} (\CTKZ) The equation $(\p-\pi(L_{Gaudin}(z)) \Psi(z)=0$
is precisely the particular case of the KZ-equation.\footnote{From the point of view of the Gaudin system
the first space $V_{0}=\CC^n$ - is  the
"auxiliary space"; and the tensor product of the other spaces $V_1\otimes ... \otimes V_\mpn$ is the
 "quantum space"
(using integrability theory language). The auxiliary and quantum spaces play the very different roles
in integrability. But in KZ-equation all the spaces $V_i,~ i=0,...$ play the same role.
This uncovers unexpected symmetry between auxiliary and quantum spaces. It is different
from the other relations between KZ-equation and the Gaudin model (e.g. \FFR).}
\\
Here $L_{Gaudin}(z) =
\sum_{i=1...\mpn} \frac{ \sum_{ab}E_{ab}\otimes \pi_i(e_{ab}^{(i)})}
{z-z_i}$ is the standard Lax matrix of the quantum Gaudin system
(see formula \ref{Lax-st} page \pageref{Lax-st}) considered in a representation $\pi$;
$(\pi, V_1\otimes ...
\otimes V_\mpn)$ is a representation of $(gl_n\oplus ... \oplus gl_n)$;
$\Psi(z)$ is a $\CC^n\otimes V_1\otimes ...\otimes V_{\mpn}$-valued
function, so that its components $\Psi_i(z)$  are $V_1\otimes
...\otimes V_{\mpn}$-valued functions; $E_{ab}\in Mat_\ms$ and $e_{ab}\in gl_{\ms}\subset
U(gl_{\ms})$ are the standard matrix units.

{\Prop
Let $\Psi(z)$ be a solution of the KZ-equation with $\kappa=1$
(\ref{KZ-fml}), then: \bea \forall i=1,...,n~~~
\pi(det(\p-L_{Gaudin}(z))) \Psi_i(z)=0 \eea } \PRF The adjoint
matrix $(\p-L_{Gaudin}(z))^{adj}$ exists as we discussed above (here
we use $\kappa=1$), so:
\begin{equation}\begin{split}
&\pi(\p -  L_{Gaudin}(z))\Psi(z)=0, \Rightarrow
\pi((\p-L_{Gaudin}(z))^{adj}) \pi(\p -  L_{Gaudin}(z))\Psi(z)=0,
\\&\text{hence }\quad
\pi( det(\p-L_{Gaudin}(z)) ) Id~~ \Psi(z)=0, \\&\text{hence  }
\forall i=1,...,n \quad \pi( det(\p-L_{Gaudin}(z))
)\Psi_i(z)=0.
\end{split}\end{equation}
\hfill\BX

{\Rem ~} The equation $det(\p -L(z)) Q(z)=0$ should be seen as
generalized Baxter's $T-Q$-equation (see \CTBig), it is well-known
that Baxter's equation plays a crucial role in the solution of
quantum systems. In particular it is important to establish that its
solutions  are rational functions (see \MTVShap). The proposition
above relates this question to the rationality of KZ-solutions
which is a more natural question (see \Sah06).

\subsection{Inversion of \MMs \label{Inv-c-ss} }
{\Th \label{Inverse-conj}
Let $M$ be a
\GT\  matrix, and assume that a two sided inverse matrix $M^{-1}$
exists
(i.e. $M^{-1} M=MM^{-1}=1$). Then $M^{-1}$ is again
a \GT\  matrix.
}
\\

\PRF The proof of the assertion consists in a small extension
and rephrasing of the arguments in the proof of Lemma 1 page 5 of
\BabelonTalon02.

We consider the Grassman algebra $\CC[\psi_1,...,\psi_n]$ (i.e.
$\psi_i^2=0, \psi_i \psi_j= -\psi_j \psi_i$), with $\psi_i$ commuting
with $M_{p,q}$. Let us introduce new
variables $\tilde \psi_i$ by: \bea (\tilde \psi_1, ... ,  \tilde
\psi_n) = (\psi_1, ... ,   \psi_n) \left(\begin{array}{ccc}
 M_{1,1} & ... &  M_{1,n} \\
... \\
 M_{n,1} & ... &  M_{n,n}
\end{array}\right).
\eea

It is easy to see that the Manin relations $[ M^{-1}_{ij}
M^{-1}_{kl}]=[M^{-1}_{kj}, M^{-1}_{il}]$ we have to prove follow from
the equality
\bea
\label{LemInvReform-fml}
 M^{-1}_{ij} M^{-1}_{kl} \psi_1 \wedge ... \wedge \psi_n=
M^{-1}_{kj} M^{-1}_{il} \psi_1 \wedge ... \wedge \psi_n+ \nn\\ +
(det(M))^{-1} \tpsi_1 \wedge \tpsi_2 \wedge ...
\wedge \stackrel{i-th~place}{\psi_j} \wedge ...\wedge \stackrel{k-th~place}{\psi_l} \wedge ...\wedge \tpsi_n.
\eea (Here assume $i \ne k$; otherwise the desired result is a
tautology). Let us prove (\ref{LemInvReform-fml}). By the definition
of the $\tpsi_i$'s, $ (\tilde \psi_1, ... ,  \tilde \psi_m) =
(\psi_1, ... ,   \psi_n) M$, multiplying this relation by $M^{-1}$
on the right we get
\[
(\tilde \psi_1, ... ,  \tilde \psi_n)M^{-1} = (\psi_1, ... , \psi_n)
,\quad \Leftrightarrow\quad  \sum_v \tpsi_v M^{-1}_{vl} = \psi_l.\] We
multiply this relation with the $n-1$-vector $\tpsi_1 \wedge \tpsi_2
\wedge ...\wedge \stackrel{i-th~place}{\psi_j} \wedge ...\wedge
\stackrel{k-th~place}{ empty } \wedge ...\wedge \tpsi_n$, and use
the fact that $\forall m: \tpsi_m^2=0$, to get
\[
( \tpsi_i M^{-1}_{il} +  \tpsi_k M^{-1}_{kl}  - \psi_l) \tpsi_1
\wedge \tpsi_2 \wedge ...\wedge \stackrel{i-th~place}{\psi_j} \wedge
...\wedge \stackrel{k-th~place}{ empty } \wedge ...\wedge \tpsi_n
=0.\] Now, using:  $det(M) M^{-1}_{ji} ~ \psi_1 \wedge ... \wedge
\psi_n = \tpsi_1 \wedge ...  \wedge \stackrel{j-th~place}{\psi_i}
\wedge ...   \wedge \tpsi_n$ we get the desired result.\BX {\Rem
~\label{konv-rem-1}}
The paper \KonvalinkaA~ 
(see section 5 page 11 proposition 5.1) contains a somewhat weaker
proposition for matrices of the form $1- tM$, where $t$ is a formal
parameter. It states that  for $(1- tM)^{-1}$ the cross-commutation
relations holds true. The "column commutation" property is not
considered there.
The equality  $det^{column} ((1-tM)^{-1})=(det^{column}(1-tM))^{-1}$
is contained in theorem 5.2 page 13 of \KonvalinkaA. The  proofs
presented there are based on deep combinatorial arguments.

\subsubsection{\label{ERBT-ss} Application to
the Enriquez - Rubtsov - Babelon - Talon theorem.} Let us show that
the remarkable theorem \EnriquezRubtsov01 (theorem 1.1 page 2),
\BabelonTalon02 (theorem 2 page 4) about "quantization" of
separation relations -- follows as a particular case from theorem
\ref{Inverse-conj} above.

Let us briefly recall the constructions of this Theorem (following,
for simplicity \BabelonTalon02). Let $\{\alpha_i,\beta_i\}_{i=1,\ldots, g}$ be
a set of quantum
``separated'' variables, i.e. satisfying the commutation relations
\[
[\alpha_i, \alpha_j ] = 0, [\beta_i, \beta_j ] = 0, [\alpha_i, \beta_j
] = f(\alpha_i, \beta_i)\delta_{ij},~~~ i,j,=1\ldots,g.
\]
satisfying a set of 
equations (i.e., quantum Jacobi
separation relations) of the form
\begin{equation}\label{BTBe}
\sum_{j=1}^g B_j(\alpha_i, \beta_i)H_j + B_0(\alpha_i, \beta_i) = 0,
i=1,\ldots, g,
\end{equation}
for a suitable set of quantum Hamiltonians $H_1,\ldots,H_n$. One
assumes that some ordering in the expressions $B_a, a=0,\ldots g$
between $\alpha_i,\beta_i$ has been chosen, and that the operators $H_i$
are, as it is written above, on the right of the $R_j$. Also, the
$B_a$ are some functions (say, polynomials) of their arguments with
$\CC$-number coefficients. Then the statement is that
the quantum operators $H_1,\ldots,H_n$  fulfilling
(\ref{BTBe}) commute among themselves).

Our proof starts from the fact that one can see that 
equation (\ref{BTBe}) can be compactly written, in matrix form, as
\[
B\cdot H=-V, \quad \text{with } B_{ij}=R_j(\alpha_i,\beta_i),
V_i=B_0(\alpha_i,\beta_i),
\]
and thus  (making contact with the formalism of
\EnriquezRubtsov01) one is lead to consider the $g\times (g+1)$ matrix
\[
A=\left(\begin{array}{cccccccc}V_1 & B_{1,1} & ... & B_{1,g} \\
 ... & ... &... & ...  \\
V_{g} & B_{g,1} & ... &  B_{g,g}
\end{array}\right)
\]
Thanks to the functional form of the $B_{ij}$'s and of the $V_i$'s,
this matrix is a Cartier-Foata matrix (i.e., elements form different
rows commute among each other), and hence, {\em a fortiori} a \MM.

Given such a $g\times (g+1)$ Cartier-Foata matrix $A$, we consider
the following  $(g+1)\times (g+1)$ Cartier-Foata (and hence, Manin)
matrix:
\[
 \tilde A=\left(\begin{array}{cccccccc}
1 & 0 & ... & 0 \\
V_1 & B_{1,1} & ... & B_{1,g} \\
 ... & ... &... & ...  \\
V_{g} & B_{g,1} & ... &  B_{g,g}
\end{array}\right)
\]
Now it is obvious that the solutions $H_i$ of 
equation
(\ref{BTBe}) are the elements of the first column of the inverse of
$\tilde A$, and namely
\[
H_i=(\tilde{A}^{-1})_{i+1,1}, i=1,\ldots,g.
\]
Since $\tilde A$ is Cartier-Foata, its inverse is Manin, and thus
the commutation of the $H_i$'s can be obtained from theorem
\ref{Inverse-conj}.\hfill\hskip-20pt\BX.

{\Rem ~} For the sake of simplicity, we considered, as in \BabelonTalon02, an
arbitrary system for which a quantized spectral curve exists, and
the commuting Hamiltonians can be written by the formula above in
terms of the so-called separated variables (see section \ref{SV-s}).
However, it should be noticed that the theorem hold in a more
general (in classical terms, St\"ackel or Jacobi) setting, where the
functional dependence of the matrix elements $B_{ij}(\alpha_i,\beta_i)$
and $V_i(\alpha_i,\beta_i)$ on their argument may depend on the index
$i$.

\subsection{Block matrices and Schur's complement\label{Schur-ss}}
{\Th \label{det-block-prop1}
Consider a 
\MM~ $M$  of size $n$, and denote its block as follows:
\bea \label{Schur-conj-fml1}
M=\left(\begin{array}{cc}
A_{k\times k} & B_{k\times n-k}\\
C_{n-k\times k} & D_{n-k\times n-k} \\
\end{array}\right)
\eea
Assume that $M,A,D$ are invertible i.e.  $\exists M^{-1},A^{-1},
D^{-1}:$ $A^{-1}A=AA^{-1}=1$, $D^{-1}D=D D^{-1}=1$, $M^{-1}M=M
M^{-1}=1$. Then the same formulas as in the commutative case hold,
namely: \begin{equation} \label{block-det} det^{col}(M)= det^{col}(A)
det^{col}(D- C A^{-1} B)
=det^{col}(D) det^{col}(A- B D^{-1} C), \end{equation}
and, moreover, the Schur's complements:  $D- C A^{-1} B$, and $A-
B D^{-1} C$  are \MMs.}

\SPRF The full proof of this fact will be given in \cite{CF}.
Actually it is not difficult to see
that $A- B D^{-1} C$, $D- C A^{-1} B $ are \MMs~ from theorem
\ref{Inverse-conj}, so basically we need to prove the determinantal
formula
(\ref{block-det}).
To this end, one proves the
{\Lem \label{block-mult-prop1}
 Let $M$ be a \MM~ of size $n\times n$. Let $X$ be a $k\times (n-k)$ matrix, for some $k$,
with arbitrary matrix elements\footnote{No conditions of commutativity at all
  are required.}.
{\bf Then:}
\bea
det^{column} ( M
\left(\begin{array}{cc}
1_{k\times k} &  X_{k\times n-k}\\
0_{n-k\times k}  & 1_{n-k \times n-k} \\
\end{array}\right) ) =det^{column} M.
\eea
}
Now, the desired result follows now from the equality:
\bea
\left(\begin{array}{cc}
A & B\\
C & D \\
\end{array}\right)
\left(\begin{array}{cc}
1 & -A^{-1} B\\
0 &  1 \\
\end{array}\right)
=
\left(\begin{array}{cc}
A & 0\\
C & D-CA^{-1} B \\
\end{array}\right).
\eea

{\Rem ~}
1) Alternatively, one can prove (\ref{block-det}) using
the Jacobi type theorem by I. Gelfand, V. Retakh
(see e.g. \GR97 theorem 1.3.3, page 8)
with similar arguments to those used
by D. Krob, B. Leclerc \KL94 (theorem 3.2 page 17).
\\
2)
The Jacobi ratio theorem was
proved for \MMs~ of the form $1- tM$, $t$ is a formal parameter,
in the remarkable paper \KonvalinkaA
(see theorem 5.2 page 13). It is actually equivalent to our
formula \ref{block-det} for matrices of such a form.
As in the case of the matrix inversion formula,
column commutativity is not considered, and  the proofs therein
contained are based on combinatorial properties.

\subsubsection{Application. MTV-Capelli identity for $gl[t]$ (for the
Gaudin system) \label{Cap-sssect} }
In \MTVCap~ a remarkable generalization of the Capelli identity has been found.
Using the beautiful insights contained therein, we remark that it
follows from theorem \ref{det-block-prop1}.

Consider $\CC[p_{i,j}, q_{i,j}]$,  $i=1,...,\ms;j=1,...,\mpn$,
endowed with the standard Poisson bracket:

$\{  p_{i,j},  q_{k,l} \}=\delta_{i,k} \delta_{j,l}$,
$\{   p_{i,j},   p_{k,l} \}= \{  q_{i,j},  q_{k,l} \}=0$.
Consider its quantization i.e. the standard Heisenberg (Weyl) associative algebra generated by
$\hat p_{i,j},\hat q_{i,j}$,  $i=1,...,\ms;j=1,...,\mpn$,
with the relations
$[ \hat p_{i,j}, \hat q_{k,l} ]=\delta_{i,k} \delta_{j,l}$,
$[ \hat  p_{i,j}, \hat  p_{k,l} ]= [ \hat q_{i,j},\hat  q_{k,l} ]=0$.
Define $n\times k$-rectangular matrices $Q_{classical}, P_{classical}, \hat Q,\hat P$
as follows with: $ Q_{i,j}= q_{i,j}$, $\hat Q_{i,j}=\hat q_{i,j}$, $P_{i,j}= p_{i,j}$,
$\hat P_{i,j}=\hat p_{i,j}$.
Let $K_1, K_2$  be  $n\times n, k \times k $ matrices with elements in $\CC$.

Let us introduce the following notations:
\bea \label{Lax4Cap}
L^{quantum}(z)= K_1+ \hat Q (z-K_2)^{-1} \hat P^t,\>
 L^{classical}(z)= K_1+ Q_{classical} (z-K_2)^{-1} P_{classical}^t
\eea
These are  Lax matrices of Gaudin type i.e. they satisfy the commutation
relation
(\ref{Lax-G-com-rel}).

Let us reformulate the result of \MTVCap~
in the following simple form:
{\Prop Capelli identity for the $gl[t]$ (Gaudin's system) }
\bea \label{Cap-Wick-1}
Wick(det(\l-L^{classical}(z)))=det(\p-L^{quantum}(z)) \eea Here we
denote by $Wick$ the linear map: $\CC[\l, p_{i,j}, q_{i,j}](z) \to
\CC[ \p, \hat p_{i,j},\hat q_{i,j}](z)$, defined as: \bea
\label{Wick-def} Wick(  f(z) \l^a \prod_{ij} q^{c_{ij}} \prod_{ij}
p_{ij}^{b_{ij}}   ) = f(z) \p^a \prod_{ij} \hat
q_{ij}^{b_{ij}}\prod_{ij} \hat p_{ij}^{b_{ij}} \eea i.e. that, from
a commutative monomial, makes a noncommutative monomial according to
the rule that any $\hat q $ placed on the left of any $\hat p$, and
the same for $z$ and $\p$. It is sometimes called "Wick or normal
ordering" in physics.
\\
\PRF
Following \MTVCap\
we consider the following block matrix: \bea
MTV=\left(\begin{array}{cc}
z  -K_2 & \hat P^t\\
\hat Q & \p-K_1 \\
\end{array}\right)
\eea It is easy to see that $MTV$ is a \MM (indeed it is strictly
related to the Manin matrices briefly introduced in (\ref{MTVex})) ;
also, as it was observed by Mukhin, Tarasov and Varchenko, the Lax matrix of
the form
(\ref{Lax4Cap}) appears as the Schur's complement $''D-CA^{-1}B''$ of the
matrix MTV:
\[
\p-K_1- \hat Q (z-K_2)^{-1} \hat P^t =\p
-L^{quantum}(z).\]
By theorem \ref{det-block-prop1} we get that
\begin{equation}
\label{Cap-MTV-1} det^{column}(MTV)=det(z-K_2) det^{column}(\p-K_2-
\hat Q (z-K_2)^{-1} \hat P^t). \end{equation} This a form of the
Capelli identity as explained in \MTVCap.

In order to arrive at formula (\ref{Cap-Wick-1}) above, we just
remark that in $det^{column}(MTV)$ all  variables $z, \hat q_{ij}$
stands on the left of the variables $\p, \hat p_{ij}$. This is due
to the column expansion of the determinant, where at first appears
the first column, then the second, and so on an d so forth.
Actually, the operators $z, \hat q_{ij}$ stand in the first $n$-th
columns of $MTV$ matrix, while the operators $\p, \hat p_{ij}$ are
in the m rightmost columns. So dividing (\ref{Cap-MTV-1}) by
$det(z-K_2)$ we obtain  (\ref{Cap-Wick-1}).
\hfill\BX
See section \ref{Cap-fml-ss} and the references quoted therein for
the relation with the classical Capelli identity.
\section{Spectral properties of \MMs\ and applications}
\subsection{The  Cayley-Hamilton theorem \label{CH-ss} }
The Cayley-Hamilton theorem, (i.e, that any ordinary matrix
satisfies its characteristic polynomial) can be considered as one of
the basic results in linear algebra. The same holds - with a
suitable proviso in mind - for \MMs.
{\Th \label{CH-th} Let $M$ be a
$n\times n$ \GT\  matrix.
Consider its characteristic polynomial and the
coefficients $h_i$ of its expansion in
powers of $t$: $det^{column}(t-M)=\sum_{i=0,...,n} h_i t^i$; then
\bea
\sum_{i=0,...,n} h_i M^i=0 \mbox{ ~i.e. ~ } det^{column}(t-M)|_{t=M}^{right~substitute}=0.
\eea
If $M^t$ is a \MM, then one should use {\em left} substitution and row determinant:\\
$det^{row}(t-M)|_{t=M}^{left~substitute}=0$.}

\PRF
It was proved in Proposition \ref{ad-inv} that there exists an
adjoint matrix $ (M-t~Id)^{adj}$, such that \bea (M-t~Id)^{adj}
(M-t~Id)=det(M-t~Id) Id. \eea The standard idea of proof is very
simple: we want to substitute $M$ at the place of $t$; the
LHS of this equality vanishes manifestly, hence we obtain the
desired equality $det(M-t~Id)|_{t=M}=0$. The only issue we need to
clarify is {\em how} to substitute $M$ into the equation and why the
substitution preserves the equality.

Let us denote by $Adj_k(M)$ the matrices defined by:
$\sum_{k=0,...,n-1} Adj_k(M) t^k= (M-t~Id)^{adj}$. The equality
above is an equality of polynomials in the variable $t$:
\bea
(\sum_k Adj_k(M) t^k) (M-t ~ Id)= \sum_k Adj_k(M) M  t^k - \sum_k Adj_k(M)   t^{k+1} =
\\ = det (M-t~Id)= \sum_k h_k t^k.
\eea
This means that the coefficients
of $t^i$ of both sides of the relation coincide. Hence we can
substitute $t=M$ in the equality, substituting "from the right":
\bea \sum_k Adj_k M  M^k - \sum_k Adj_k M^{k+1} = \sum_k h_k
t^k\big\vert_{t=M}. \eea The left hand side is manifestly zero, so
we obtain the desired equality:
$\sum_{k=0,...,n} h_k M^k=0.$\hfill\BX.
{\Rem ~} See  \Molev02 (section 4.2), \GIOPSSupB~ and references therein
for other CH-theorems.

\subsubsection{Example. The Cayley-Hamilton theorem for the Yangian
 \label{Yang-CH}
 }
Let us consider the Lax matrix $T_{gl_n-Yangian}(z)$ (or $T(z)$ for brevity)
of Yangian type (see section \ref{Yang-ss}).
The matrix  $T(z)e^{-\p}$ is a \MM. Let us derive a Cayley-Hamilton identity
for $T(z)$ from the one for the \MM\ $T(z)e^{-\p}$.

{\Def Let us define the "quantum powers"
$T^{[p]}(z)$ for the Yangian type Lax matrix as follows:
\bea \label{QP-Yang-fml}
T^{[p]}(z)\stackrel{def}{=} T(z+p)\cdot T(z+p-1)\cdot\ldots \cdot
T(z+1)=e^{p\p} (T(z)e^{-\p})^p
\eea
}
Denote by $QH_i(z)$ \footnote{we slightly changed the definition of $QH_i(z)$ comparing with the version 1,
and correct the misprint in formula \ref{f001}}
coefficients of the expansion of the ``ordered''
characteristic polynomial in powers of $e^{\p}$:
\bea \label{f001}
det^{column}(t-T(z)e^{-\p})=\sum_{k=0,...,n} t^k  QH_k(z) e^{(k-n)\p},
\eea
Note that the operator $\p$ does not enter into the expressions
$QH_i(z)$.
Explicit formulas for $QH_i(z)$ are obviously the following:
\bea
& QH_n(z) = 1, ~~~  QH_{n-1} (z) =- Tr(T(z)), ~~~ QH_0(z)= (-1)^n~ qdet(T(z)), & \\
&\qquad QH_{n-i}(z)= (-1)^{(n-i)}~
\sum_{j_1<...<j_i} qdet( T(z)_{j_1,....,j_i} ).&
\eea
In words, the
$QH_i(z)$'s are the sums of principal q-minors of size $i$ (in
complete analogy with the commutative case, modulo substitution of
minors by q-minors). Recall that  $qdet(M(z))$
is defined by the formula $qdet(M(z))=\sum_{\sigma\in S_n} \prod_k
M_{\sigma(k),k}(z-k+1)$ (see e.g. formula 2.24 page 10 \Molev02).

{\Prop \label{prop-Yang-CH}
 The following is an analogue of the Cayley-Hamilton theorem
for Yangian type Lax matrices: \bea \sum_{k=0,...,n} QH_k(z+n)
T^{[k]}(z)=0. \eea }
\PRF The proof is very simple: one uses
Cayley-Hamilton theorem  \ref{CH-th} for the \MM~ $T(z)e^{-\p}$, and
than excludes powers of $e^{\p}$ from the identity\footnote{One can
easily exclude $\p$ in the Yangian case, but it is not clear at the
moment how to do it in the Gaudin case. Nevertheless, a similar
identity can be proved in for the Gaudin case also: \CTKZ~ (by
different methods).}.
We can proceed as follows. By theorem
\ref{CH-th} it is true:
\bea \label{id0} \sum_{k=0,...,n} h_k(z,\p)
(T(z)e^{-\p})^{k}=0, ~~~~~ \Rightarrow ~~~~~ \sum_{k=0,...,n}
h_k(z,\p) e^{-k\p}T^{[k]}(z)=0,
\eea
where $h_k(z,\p)$ are defined
as $det^{column}(t-T(z)e^{-\p})=\sum_{i=0,...,n} h_i(z,\p) t^i$.
By the definition of $QH_i(z)$ (formula \ref{QP-Yang-fml}) one has: $h_i(z,\p)=QH_i(z) e^{(i-n)\p}$,
substituting this into \ref{id0}, one obtains:
\bea \label{id00}
0=\sum_{k=0,...,n} h_k(z,\p) e^{-k\p}T^{[k]}(z)
= \sum_{k=0,...,n}  QH_k(z) e^{(k-n)\p} e^{-k\p}T^{[k]}(z) = \\
= \sum_{k=0,...,n}  QH_k(z) e^{-n\p} T^{[k]}(z)
=  e^{-n\p} \sum_{k=0,...,n}   QH_k(z+n)  T^{[k]}(z)
,
\eea
So we conclude that $ \sum_{k=0,...,n}   QH_k(z+n)  T^{[k]}(z) =0$.
\hfill\hskip-18pt\BX
{\Rem ~} The quantities $QH_k(z)$ were defined in %
\cite{KulishSklyanin81} (formula 5.6 page 114) by means of different
arguments. The idea of obtaining them as coefficients of a suitable
characteristic polynomial is due to D.Talalaev \T04.

\subsection{The Newton identities and applications \label{Newton-ss}}
As it is well known, the Newton identities
are identities between power sums  symmetric functions and
elementary symmetric functions. They can be rephrased as relations
between $TrM^k$ and coefficients of $det(t+M)$ for matrices with
commutative entries. We will show that the same identities hold true
for \MMs~ and  present  applications.
{\Th \label{Newt-th} Let $M$ be a 
\GT\  matrix, and denote by $ h_k$ the coefficients of the
expansion in $t$ of the determinant of $t-M$, i.e.,
$det^{column}(t-M)= h_k t^k$; conventionally, let $h_k=0$ for $k<0$.
Then the following identity holds:
\bea
 \partial_t det^{column}(t-M) =\frac {1}{t} (det^{column}(t-M))
 \sum_{k=0,...,\infty} Tr (M/t)^k,\\
\label{Newt-fml2}
\Leftrightarrow \forall k: -\infty <  k\le n \mbox{ ~ it holds: ~ } ~~~~~
k  h_k=\sum_{i: max(0,-k)\le i \le n-k }  h_{k+i} Tr (M)^{i},
\eea
where $n$ is the size of the matrix $M$. If $M^t$
is a \MM, then $\partial_t det^{row}(t-M) = (det^{row} (t-M))
 \left( 1/t \sum_{k=0,...,\infty} Tr (M/t)^k \right) $.
}

So these identities are identical to those of the commutative case,
provided one pays attention to the order of terms: $\tilde h_i Tr
M^p$ if $M$ is a \MM\
($Tr M^p \tilde h_i $ if $M^t$ is a \MM).
Obviously enough, this difference is due to the absence of the
commutativity.

\SPRF First observe the following simple property, whose proof is
straightforward\footnote{Actually, the equality holds for any matrix
$M$, provided one consistently defines the determinant and the
adjoint matrix.}:
\[
Tr (t+M)^{adj}=\partial_t det(t+M).
\]
So we consider
\bea 1/t \sum_{k=0,...,\infty} Tr\big((-M/t)^k\big) = Tr
\frac{1}{t+M}=Tr \Bigl((det^{column}(t+M))^{-1}(t+M)^{adj} \Bigr) =
\nn\\= (det^{column}(t+M))^{-1} Tr (t+M)^{adj}=
(det^{column}(t+M))^{-1}
\partial_t det^{column}(t+M).
\eea
Substituting $-M$ instead of $M$ one obtains the result.
The case where
$M^t$ is a \MM~ is similar. \hfill \BX.
{\Rem ~} While working at the proof of this proposition we were
informed by P. Pyatov (unpublished, based on \GIOPS) that the Newton
identities holds true as well as their q-analogs. Our proof is
different and very  simple, but it is a challenge how to generalize
such an argument to the case  $q\ne 1$.

\subsubsection{Quantum powers,
$ZU_{crit}(gl_n[t,t^{-1}],c)$, and new quantum Gaudin
integrals\label{QP-G-sss}}
We recall from
\CTKZ ~
a possible definition of "quantum powers" of the Lax matrix of the Gaudin
model and
apply the  quantum Newton relations above to prove the commutativity
of their traces. Hence we obtain new explicit generators in the
center of $U_{c=crit}(gl_n[t,t^{-1}],c)$ and  in the commutative
Bethe subalgebra in $U(gl_n[t]))$ (see Appendix (section
\ref{App-s}) for a reminder); in a physicist's language we give
explicit formulas for quantum conservation laws for the quantum
Gaudin system. They are quantum analogs of $Tr (L^{classical}
(z))^k$, and the result below exhibits the appropriate corrections
for traces of powers which preserves commutativity.

{\Def The quantum powers of Gaudin type
Lax matrices (definition \ref{GL-def},
formula \ref{Lax-G-com-rel})
are defined inductively as follows:
\vskip -5mm
\bea
\label{quant-power}
&&~~~~~~~~~~~~~~~~~~~~L^{[0]}(z)=Id, ~~~~
L^{[i]}(z)= L^{[i-1]}(z)L(z)+ (L^{[i-1]}(z))'
\eea } Here $(L^{[i-1]}(z))'$ is the derivative with respect to $z$ of
$(L^{[i-1]}(z))$. $L^{[i]}(z)$ are noncommutative analogues of the
Fa\`a di Bruno polynomials (\cite{Di03} section 6A, page 111).
(Remark that in the commutative case $L(z)$ and $L(z)'$ commute).
The binomial type formula below is due to D. Talalaev:
\bea
(\p-L(z))^i= \sum_{p=0,...,i}
(-1)^p \left(\begin{array}{c} i \\ p
\end{array}\right)
  \p^{(i-p)} L^{[p]}(z).
\eea
{\Th ~ \label{NewtProp1}  Consider the quantum Hamiltonians
$QH_i(z)$ defined by:
$det(\p-L(z))=\sum_{i=0,...,n} QH_{i}(z) \partial_z^i$. Then:
\bea
\forall k,l=1,\ldots,n\text{ and } u,v\in \CC ~~~~
[QH_k(u) , Tr L^{[l]}(v) ]=0, \quad
[Tr L^{[k]}(u), Tr L^{[l]}(v) ]=0.
\eea
}

\SPRF
The proof follows immediately from the  Newton identities above,
theorem of \T04 ~
(see formula \ref{TT}) on the
commutativity of coefficients of $det(\p-L(z))=\sum_i  QH_i(z)\p^i$,
and the binomial formula above. \hfill\BX
{\Cor \label{Tr-Bet} If $L(z)=L_{gl_n[t]}(z)$ is given by formula
(\ref{Lax-glt}b), then $Tr L^{[k]}_{gl_n[t]}(z)$ are generating
functions (in the variable $z$) for the elements in a commutative
Bethe subalgebra in $U(gl_n[t])$. Such elements for $k=1,...,n$
algebraically generate the Bethe subalgebra. }
{\Cor If \label{Tr-cen} $L_{full}(z)$ is given by formula
(\ref{L-full}), then $:Tr L^{[k]}_{full}(z):$ are generating
functions (in the variable $z$) for the elements in the center of
$U_{c=crit}(gl_n[t,t^{-1}],c)$. Such elements for $k=1,...,n$
algebraically generate the center. } (See Appendix (section
\ref{App-s}) for the definitions of $U_{c=crit}(gl_n[t,t^{-1}],c)$,
the Bethe subalgebra and the normal ordering of operators
$:\cdots:$.)
\subsubsection{Application. The Newton identities for the Yangian
\label{QP-Y-sss}
}
 We recall that if a  matrix $T(z)$ is a $n\times
n$ Lax matrix of Yangian type
(see section \ref{Yang-ss}),%
then the matrix  $T(z)e^{-\p}$ is a \MM.
We have defined (see formula \ref{QP-Yang-fml})
quantum powers for the Yangian type Lax operators as follows:
\bea
T^{[p]}(z)\stackrel{def}{=} T(z+p)T(z+p-2)\times ... \times T(z+1)=e^{p\p} (T(z)e^{-\p})^p.
\eea
Writing the Newton identities for $T^{[p]}(z)$ (see formula
\ref{QP-Yang-fml}) we obtain generators in the
commutative subalgebra of quantum integrals of motion that, to the
best of our knowledge, were not known in the literature.

With the notations of section \ref{Yang-CH} page \pageref{Yang-CH}
we denote by $QH_i(z)$ coefficients of the expansion of the
``ordered'' characteristic polynomial in powers of $e^{\p}$: \bea
\label{f001b} det^{column}(t-T(z)e^{-\p})=\sum_{k=0,...,n} t^k
QH_k(z) e^{(k-n)\p}, \eea Note that the operator $\p$ does not enter
into the expressions
$QH_i(z)$.
Explicit formulas for $QH_i(z)$ are
\[QH_0(z)= (-1)^n~ qdet(T(z)),\quad
QH_{n-i}(z)= (-1)^{(n-i)}~ \sum_{j_1<...<j_i} qdet( T(z)_{j_1,....,j_i} ).\]
 Recall that  $qdet(M(z))$
is defined by the formula $qdet(M(z))=\sum_{\sigma\in S_n} \prod_k
M_{\sigma(k),k}(z-k+1)$, and that the
$QH_i(z)$'s are the sums of principal q-minors of size $i$ (in
complete analogy with the commutative case, modulo substitution of
minors by q-minors).

{\Th ~  \label{Newt-Yang-th}
\vskip -12mm
\bea \label{qTr-com}
[   QH_i(z) , Tr T^{[i]}(u) ]=0, ~~~~~~ [ Tr T^{[k]}(z), Tr T^{[i]}(u)]=0,
\eea
\bea \label{qTr-Newt-yang}
k  QH_k(z-k+n)  = \sum_{i: max(0,-k)\le i \le n-k } QH_{k+i}(z-k+n) Tr  T^{[i]}(z)
\eea
which is the Newton identity for the Yangian.}

\PRF The commutativity of the $  QH_i(z) $ is well-known
\cite{KulishSklyanin81} page 114 formulas 5.6-5.7,
\NO.
Hence (\ref{qTr-com}) follows from  (\ref{qTr-Newt-yang}), which, in
its turn, follows from theorem \ref{Newt-th}.
Indeed, \ref{Newt-fml2} reads:
\bea
k  h_k=\sum_{i: max(0,-k)\le i \le n-k }  h_{k+i} Tr (M)^{i}, \\
k QH_k(z) e^{(k-n)\p} = \sum_{i: max(0,-k)\le i \le n-k } QH_{k+i}(z) e^{(k+i-n)\p}
e^{-i\p} Tr  T^{[i]}(z) ,\\
k e^{(k-n)\p} QH_k(z-k+n)  = \sum_{i: max(0,-k)\le i \le n-k } QH_{k+i}(z) e^{(k-n)\p}
 Tr  T^{[i]}(z) ,\\
k e^{(k-n)\p} QH_k(z-k+n)  = e^{(k-n)\p} \sum_{i: max(0,-k)\le i \le n-k } QH_{k+i}(z-k+n)
 Tr  T^{[i]}(z).
\eea
So we conclude that: $k  QH_k(z-k+n)  = \sum_{i: max(0,-k)\le i \le n-k } QH_{k+i}(z-k+n) Tr  T^{[i]}(z)$.
 \BX
{\Rem ~} Let us note that for $k=0$ this follows from the Cayley-Hamilton theorem for the Yangian
\ref{prop-Yang-CH} page \pageref{prop-Yang-CH} by taking the trace.

\subsection{The (Quantum) MacMahon-Wronski formula}\label{MW-s}
The following identity for a $s\times s$-matrix $M$  over a commutative ring
is called the MacMahon (and sometimes Wronski (\UmedaW,
\GIOPSW ~page 9
)) formula:
\bea
1/det(1-M)=\sum_{n=0,...,\infty} Tr {S^n M},
\eea
where $S^n M$ is $n$-th symmetric power of $M$.
It can be easily verified by diagonalizing the matrix $M$.
\\
{\em {\bf Theorem}  \GLZ 
 \footnote{see the papers\GIOPS, \UmedaW, 
\EtingofPak06,\HaiLorenz06,\KonvalinkaPak06,\FoataHan06
~
for related results}
 the identity above holds true for \MMs~ (and more generally for their q-analogs)
with the following definition of $Tr{S^n M}$:} \bea Tr S^n M =
1/n!\sum_{l_1,...l_n: 1\le l_i \le n} perm^{row}
\left(\begin{array}{cccc}
M_{l_1,l_1} & M_{l_1,l_2}& ... & M_{l_1,l_n} \\
M_{l_2,l_1} & M_{l_2,l_2}& ... & M_{l_2,l_n} \\
... & ... & ... & ... \\
M_{l_n,l_1} & M_{l_n,l_2}& ... & M_{l_n,l_n}
\end{array}\right).
\eea
\noindent We remark that in this formula repeated indexes are
allowed.
We recall that the permanent of a
\MM~ was defined by  formula \ref{M-perm-f} as follows: \bea perm
M=perm^{row} M=\sum_{\sigma\in S_n} \prod_{i=1,...,n}
M_{i,\sigma(i)}. \eea This definition of traces of symmetric powers
is the same as in the commutative case, with the proviso in mind to
use row permanents.

\subsubsection{Application. $Tr~S^K L(z)$, $Z U_{crit}(gl_n[t,t^{-1}],c)$,
new Gaudin's integrals
\label{Nazar}} In this section we apply the quantum MacMahon-Wronski
relation to construct further explicit generators in the center of
$U_{c=crit}(gl_n[t,t^{-1}],c)$ and  in the commutative Bethe
subalgebra in $U(gl_n[t]))$, i.e.,  we give explicit formulas for
another set of quantum conservation laws for the quantum Gaudin
system. They are quantum analogs of $Tr S^n L^{classical} (z)$. The
elements in the center of $U(gl_n)$
introduced implicitly in \Nazarov91 ~
(the last formula - page 131, see also \UmedaW) are particular
cases of this construction.

{\Th Let $L(z)$ be a Lax matrix of Gaudin type - see definition
\ref{GL-def}. Let us define the elements $S_n(z)$ as follows: write
$Tr S^n(\p-L(z))=\sum_{k=0,...,n} c_{k,n}(z)\p^k$ and set, by
definition, $S_n(z)=(-1)^n c_{0,n}(z)$\footnote{Namely, the quantities
$S_n(z)$ are the
$0$-th order parts of the differential operators $Tr S^n(\p-L(z))$.}
. These operators commute among themselves and with the coefficients of
the quantum characteristic polynomial $QH_i(z)$, defined via
($det(\p-L(z))=\sum_i QH_i(z) \p^i$). Explicitly \bea \forall k,p\text{ and }
z,u\in \CC:
[S_k(z), S_p(u)] = 0, ~~~ [S_k(z), QH_p(u)] = 0. \eea Hence they
provide new generators among quantum commuting integrals of motion
of the Gaudin system. } {\Cor \label{Naz-Bet} Consider
$L(z)=L_{gl_n[t]}(z)$ given by formula \ref{Lax-glt}b. Then $S_n(z)$
are generating functions in the variable $z$ for the elements in the
commutative Bethe subalgebra in $U(gl_n[t])$. Such elements for
$k=1,...,n$ algebraically generate the Bethe subalgebra. }

{\Cor Consider \label{Naz-cen} $L_{full}(z)$ given by formula \ref{L-full};
then $:S_k(z):$ are generating functions in the
variable $z$ for the elements in the center
of $U_{c=crit}(gl_n[t,t^{-1}],c)$.
Such elements for $k=1,...,n$ algebraically generate the center.
}
(See Appendix (section \ref{App-s}) for definitions
of $U_{c=crit}(gl_n[t,t^{-1}],c)$, the Bethe subalgebra and the
normal ordering $:...:$.)

\SPRF
The proof follows from the MacMahon-Wronski relations and
a theorem of \T04 ~
(see formula \ref{TT}) on commutativity of coefficients of $det(\p-L(z))$.
Corollary \ref{Naz-Bet} follows immediately, corollary \ref{Naz-Bet} with the help
of theorem \ref{Z-th-app} (see Appendix (section \ref{App-s})).
\hfill\BX
{\Rem ~}
Applying the same construction to the Yangian case i.e. $Tr S^n (  T(z) e^{-\p})$
one obtains expressions commuting among themselves and with coefficients of $det(1-  T(z)e^{-\p} )$.
Actually these expressions are quite well-known.

\section{Quantum separation of variables \label{SV-s}}
Let us briefly discuss some results and conjectures about the
problem of separation of variables.
E.Sklyanin (see surveys \SklA,\SklB) ~%
proposed an approach which potentially should give the most powerful
way to solve integrable systems. This program is far from being
completed.
The ultimate goal in this framework is to construct coordinates
$\alpha_i,\beta_i$ such that a joint eigenfunction of all
Hamiltonians will be presented as product of functions of one
variable:
\[
 \Psi(\beta_1,\beta_2,...) = \prod_i
\Psi^{1-particle}(\beta_i).
\]
We consider this construction at the  quantum level, trying to
extend in this realm some ideas of Sklyanin and
\AdamsHarnadHurtubise93, \cite{DD94},
\Gekhtman and others.
We have checked the validity of the conjectures below for $2\times
2$ and  $3\times 3$ cases (in particular comparing with the results
of \Sklyanin92
). Let us remind that, in the classical case, the construction of
separated variables for the systems we are considering goes,
somewhat algorithmically\footnote{We are herewith
  sweeping under the rug the problem known as ``normalization of the  Baker
  Akhiezer function".}, as follows:
\begin{description}
\item[Step 1] One considers, for a $gl_n$ model,
along the Lax matrix $L(z)$, the matrix
  $M=\lambda-L(z)$ and its classical adjoint $M^\vee$.
\item[Step 2] One takes a vector $\psi$ by means
suitable linear combination of columns (or rows) of
  $M^\vee$;
in the simplest case, one can take $\psi$ to be one of the columns,
say the last of $M^\vee$.
One seeks for pairs $(\lambda_i,z_i)$ that solve
\[
\psi_i=0, i=1,\ldots,n.
\]
\item[Step 3] To actually solve this problem, one proceeds as follows.
As each component $\psi_i$ of $\psi$ is a polynomial
of degree at most $n-1$, one can form, out of $\psi$ the
matrix $M_\psi$, collecting the coefficients of the expansion of $\psi_i$ in
powers of $\lambda$, i.e.:
\[
[M_\psi]_{j,i}=\text{res}_{\lambda=0} \psi_i \lambda^{n-j-1},\quad i=1,\ldots, n,\>
j=1,\ldots, n
\]
\item[Step 4] The separation coordinates are given by pairs $(\lambda_i, z_i)$
  where $z_i$'s are roots of $\text{Det}(M_\psi)$ and $\lambda_i$ are the
  corresponding values of $\lambda_i$
\footnote{In the quadratic $R$-matrix case, actually one
    has to take as canonical momenta, the logarithms of these $\lambda_i$.},
that
  can be obtained, e.g., via the Cramer's rule. By construction , the Jacobi
  separation relations are the equation(s) of the spectral curve, $\text{Det}(\lambda-L(z)=0$.
\end{description}
\subsection{The Yangian case}
Let $T(z)$ be a Lax matrix of the Yangian type (see section \ref{Yang-ss}), so
$(1-e^{-\p}T(z))$\footnote{in previous sections we worked with $(1-T(z)e^{-\p})$,
which is actually more well suited to the standard Yangian conventions, but let
us switch here to $(1-e^{-\p}T(z))$}
 is a \MM~and adjoint matrix can be calculated by standard formulas
(see section \ref{Cramer-s}). Consider
\bea
(1-e^{-\p}T(z))^{adjoint}=
\left(\begin{array}{ccc}
(1-e^{-\p}T(z))_{1,1}^{adjoint} & ... & (1-e^{-\p}T(z))_{1,n}^{adjoint} \\
... & ... & ... \\
(1-e^{-\p}T(z))_{n,1}^{adjoint} & ... & (1-e^{-\p}T(z))_{n,n}^{adjoint}
\end{array}\right).
\eea

Let us take an arbitrary column of this matrix, say the last column,
and define $M_{i,j}(z)$ as follows:
\bea
\left(\begin{array}{ccc}
 (1-e^{-\p}T(z))^{adjoint}_{1,n} \\
 ... \\
 (1-e^{-\p}T(z))^{adjoint}_{n,n}
\end{array}\right)
=
\left(\begin{array}{ccc}
\sum_{k=0...n-1} e^{-k\p} M_{1,k} (z) \\
 ... \\
\sum_{k=0...n-1} e^{-k\p} M_{n,k} (z) \\
\end{array}\right).
\eea

Let us call "M-matrix" the matrix of these coefficients, that is :
\bea\label{transp} M_{n\times n} (z) \stackrel{def}{=}
\left(\begin{array}{ccc}
M_{1,0} (z) & ... & M_{1,n-1}  (z)\\
 ... & ... &... \\
M_{n,0} (z) & ... & M_{n,n-1}  (z)\\
\end{array}\right)^t.
\eea
That is, $M_{ij}$ is a matrix collecting
the coefficients of expansion in
left powers of $e^{-\p}$ of the elements of
the last column of $(1-e^{-\p}T(z))^{adjoint}$
{\Prop For $n=2,3$, $M(z)$  is a \MM.(Remark the transposition in formula
(\ref{transp}))}.

\subsection{The cases $\mathbf{n=2,3}$.}
In the cases of low matrix rank, computations can be done explicitly. The
following arguments hold.
\begin{enumerate}
\item Define $B(z)=Det^{column}(M(z))$;
then
\bea
[ B(z) , B(u)]=0.
\eea
\item
Consider any root $\beta$ of the equation $B(u)=0$, (it belongs to an
appropriate algebraic extension of the original noncommutative algebra $R$).
Then the overdetermined system of, respectively, $2$ and $3$
equations for the single variable
$\alpha$ has a unique solution:
\bea
\left(\begin{array}{ccc}
M_{1,0}(z)|_{Substitute~left~z \to \beta} & ... & M_{1,n-1} (z)|_{Substitute~left~z \to \beta} \\
 ... & ... &... \\
 ... & ... &... \\
M_{n,0}(z)|_{Substitute~left~z  \to \beta} & ... & M_{n,n-1}(z)|_{Substitute~left~z \to \beta} \\
\end{array}\right)
\left(\begin{array}{c}
\alpha^{n-1} \\
\alpha^{n-2} \\
 ... \\
1 \\
\end{array}\right)
=
\left(\begin{array}{c}
0\\
0\\
..\\
0\\
\end{array}\right).
\eea
\item
Consider all the roots $\beta_i$ of the equations $B(u)=0$,
and the corresponding variables $\alpha_i$.
{Then:}
\begin{itemize}
\item The variables $\alpha_i,\beta_i$ satisfy the following commutation relations:
\bea
[ \alpha_i , \beta_j ] =-\alpha_i \delta_{i,j},\\
  ~ [ \alpha_i , \alpha_j ] =0, ~~~~~  [ \beta_i , \beta_j ] =0,
\eea
\item
$\alpha_i,\beta_i$ satisfy the "quantum characteristic equation":
\bea \label{char-id}
\forall i: ~~~ det(1-e^{-\p} T(z) )|_{Substitute~left~z \to \beta_i;\ e^{-\p} \to \alpha_i } =0
\eea
\item
If  $T(z)$ is generic, then  variables $\alpha_i, \beta_i$
are "quantum coordinates" i.e. all elements of the algebra $R$ can be expressed via
 $\alpha_i, \beta_i$ and the centre (Casimirs) of the algebra $R$.
\end{itemize}
\end{enumerate}
The proof of these statements can be done by direct calculations; we remark
that our formulas reproduce those of the paper \Sklyanin92.\\
{\bf Conjecture}:
It is natural to conjecture that the same hold for higher values of $n$, that
is, that the $M$ matrix is a \MM, and it provides quantum separated variable for
the $gl(n)$ Yangian case.

{\Rem ~} In  the classical case,
(and also for the Gaudin model) this solution of the
Separation of Variables problem can be explicitly
found in \AdamsHarnadHurtubise93.
The Poisson version of
the conjecture above about Manin properties of Yangian systems,
was not yet, to the best of our knowledge, considered
in the literature; however, it
can possibly be deduced from the results of \cite{DD94}  and \Gekhtman.

{\Rem ~}
Another open problem is a corresponding conjecture for the Gaudin
case. The conjecture holds true for $n=2$, but seems not directly
extendable to the case $n>2$. {\Rem ~}
We end this section with a few other remarks on separated variables.
Let us consider the
noncommutative algebra $R$ is the algebra which enters the
definition of $T(z)$ (see section \ref{Yang-ss}), i.e. by definition the
matrix elements of $T(z)$ belong to $R((z))$.

One can realize any representation of this algebra as a subspace
in the space of functions
$F(\beta_i)$.
Consider the joint eigen-function $\Psi(\beta_1,\beta_2,...)$ of all quantum
commuting Hamiltonians
\bea
det(1-e^{-\p} T(z)) \Psi(\beta_1,\beta_2,...) = \sum_{k=0,...,n} e^{-k\p} h_{k}(z) \Psi(\beta_1,\beta_2,...)
\eea
Remark that
$h_k(z)$ are eigenvalues of the corresponding Hamiltonians $QH_i(z)$.

It is easy to see from \ref{char-id}, that
\bea \label{sep-eq}
&& 
\forall i:~~~~ \sum_{k=0,...,n} e^{-k\partial_{\beta_i}} h_{k}(\beta_i) \Psi(\beta_1,\beta_2,...)=0.
\eea
So the equations for functions of many variables have been separated i.e.
it satisfies the system of equations in each variable $\beta_i$ separately.

Although from the property above one might expect that
\bea
&&
\Psi(\beta_1,\beta_2,...) = \prod_i \Psi^{1-particle}(\beta_i),
\eea
where each $\Psi^{1-particle}(z)$ is a solution of the equation
$\sum_k e^{-k\partial_{z}} h_{k}(z) \Psi^{1-particle}(z)=0$,
we remark that formula \ref{sep-eq} guarantees
a weaker property, namely that  $\Psi(\beta_1, \beta_2,...)$
is a linear combination of such expressions.

\section{Concluding remarks and open problems}
We postpone further discussion to the subsequent paper(s) \cite{CF}
and restrict ourselves with some brief remarks.

We think that it is very important
to  continue further the study of the notion of "quantum spectral curve"
discovered by D. Talalaev.
In particular one may hope to describe quantum action-angle variables
with some "quantum Abel-Jacobi transform" related to
Talalaev's "quantum spectral curve".
This might also be important for the general development of noncommutative
algebraic geometry. More generally one may hope to develop the methods
used in classical integrability: Lax pairs,
dressing transformations, tau-functions,
explicit soliton and algebro-geometric solutions, etc. for quantum systems.

It would be interesting to define quantum immanants ("fused T-matrices")
for  $U_{c=crit}(gl_n[t,t^{-1}],c)$ and to transfer the
remarkable properties
from \OkOlsh ~
 ($U(gl_n)$ case) to this more general case.
In particular the task is to
describe Harish-Chandra map for $ZU_{c=crit}(gl_n[t,t^{-1}],c)$.
(What is "loop version" of shifted Schurs ?).
Moving away from the critical level towards  (quantum)
W-algebras is another task.

We have presented results for $gl_n$-Gaudin and $gl_n$-Heisenberg's
XXX systems. We plan to extend our study to Q-analogues (i.e.
Heisenberg's XXZ system) and super analogue in the next future,
it would be interesting to understand a relation with more
recent quantum group like structures (\cite{F07}).
 The
challenge is to consider other systems. For example
$\p-L_{so(n)~or~sp(n)-Gaudin}(z)$ is {\bf not} \MM. However, there
are evidences that there should exists an appropriate
$"det(\p-L_{so(n)~or~sp(n)-Gaudin}(z))"$, since the existence of
such an object is somewhat the core of the Langlands
correspondence. The practical problem is, however, to find an
explicit and viable formula for the determinant. In this respect we
notice that, despite $det^{column}(\p-L_{so(n)~Gaudin}(z))$ does not
depend on the order of columns, it does not produce commutative
elements already for $so(4)$ (A.Molev, A.C.)). More in general, in
view of the existence of a huge number of quantum integrable
systems,  it seems to be natural to address the questions herewith
studied to all of them. Namely, we point out the following general
question: what are the conditions on a matrix with a noncommutative
entries such that there is a natural construction of the determinant
and linear algebra works? How efficiently can they be used in the
theory of quantum integrable systems? We plan to come back to these
and related questions in the next future.

\section{Appendix}
\subsection{The center of $U_{c=crit}(gl_n[t,t^{-1}],c)$,  Bethe subalgebras
and normal ordering
\label{App-s} 
 }
Let us recall the main result of \CTBig ~
and background. (See also \cite{CM}).\\
\bea \label{L-full}
\mbox{Consider: ~~ }
L_{full}(z)=\sum_{i=-\infty ... \infty }
\frac{1}{z^{i+1}}
\left(\begin{array}{ccc}
e_{1,1}t^i & ... & e_{1,n}t^i \\
... &  ... & ... \\
e_{n,1}t^i & ... & e_{n,n}t^i
\end{array}\right) 
\eea
Here, one has to consider
$e_{ij}t^k$ as elements of $gl_n[t]\oplus t^{-1} gl_n^{op}[t^{-1}]$,
\footnote{by  $g^{op}$ is denoted Lie algebra with an opposite commutator:
$[g_1, g_2]_{g^{op}}\stackrel{def}{=} - [g_1, g_2]_{g}$ }
rather than elements of $gl_n[t,t^{-1}]$). In such a case,
$L_{full}(z)$ is of Gaudin type,
(see definition \ref{GL-def} section \ref{Gaud-Manin-ss}),
i.e., it satisfies the relations \ref{Lax-G-com-rel}.


Via Talalaev's construction one defines
the commutative subalgebra in $U(gl_n[t]\oplus t^{-1} gl_n^{op}[t^{-1}])$, i.e.
expressions $QH_i(z)$ defined by $det^{col}(\p-L_{full}(z))=\sum_k \p^k QH_k(z)$
are generating functions in $z$ for some elements in $U(gl_n[t]\oplus t^{-1} gl_n^{op}[t^{-1}])$
which commute by Talalaev's theorem and hence generate
some commutative subalgebra.
{\Def \label{Bethe-def} This commutative subalgebra
will be called "Bethe subalgebra" in $U(gl_n[t]\oplus t^{-1} gl_n^{op}[t^{-1}])$.}
{\Def  Let us call by  "Bethe subalgebra" in $U(gl_n[t])$ the commutative subalgebra
defined in a similar way with the help of $det^{col}(\p-L_{gl_n[t]}(z))$,
(see \ref{Lax-glt}b for $L_{gl_n[t]}$).}
This Bethe subalgebra is clearly an image of the Bethe subalgebra in
$U(gl_n[t]\oplus t^{-1} gl_n^{op}[t^{-1}])$ under the natural projection
$U(gl_n[t]\oplus t^{-1} gl_n^{op}[t^{-1}])\to U(gl_n[t])$.
{\Rem ~} The name "Bethe subalgebra" was proposed in \NO ~
for a related subalgebra in
(twisted)-Yangians.
(For Yangians (not twisted) the Bethe subalgebra (without a name)
was defined in \cite{KulishSklyanin81} formula 5.6 page 114.)
%
\paragraph{Normal ordering  }

Let us recall the standard definition of the normal ordering
(see e.g. V. Kac  \cite{Kac} formula 2.3.5 page 19,    formula 3.1.3 page 37).

{\Def Let $a(z), b(z)$ be arbitrary formal power series with values in arbitrary
associative algebra (in our case: $a(z)= (L_{full}(z))_{ij}, b(z)= (L_{full}(z))_{kl}$,
for some $i,j,k,l$, (see \ref{L-full} for $L_{full}(z)$).
The normally ordered product $: a(z) b(z) :$ is defined as follows:
\bea
: a(z) b(z) : \stackrel{def}{=} a(z)_{+} b(z) + b(z) a(z)_{-},
\eea
where (e.g. V. Kac  \cite{Kac} formula 2.3.3 page 19),
\bea
a(z)_{+}\stackrel{def}{=} \sum_{i\ge 0} a_{-i-1}z^i = \sum_{n<0} a_{n}z^{-n-1}, \\
a(z)_{-}\stackrel{def}{=} \sum_{i< 0} a_{-i-1} z^i = \sum_{n\ge 0} a_{n} z^{-n-1}.
\eea
}

{\Def (e.g. V.Kac \cite{Kac} formula 3.3.1 page 42)
The normally ordered product of more than two series  
$a^{1}(z), a^{2}(z), ... , a^{n}(z)$ is defined inductively
``from right to left'':
\bea
: a^{1}(z) a^{2}(z) ...  a^{n}(z) : \stackrel{def}{=}
: a^{1}(z) ...  : a^{n-1}(z) a^{n}(z) : ... :.
\eea
}

\paragraph{ The center of $U_{c=crit}(gl_n[t,t^{-1}],c)$ and the "critical level" }

Let $(gl_n[t,t^{-1}],c)$ be the  central extension of the Lie algebra $gl_n[t,t^{-1}]$:
\bea
[g_1 t^k , g_2 t^l]= [g_1, g_2] t^{k+l} +  c (n Tr(g_1 g_2)  - Tr(g_1)Tr(g_2) ) k\delta_{k,-l}.
\eea
Note that for $sl_n$ the term $Tr(g_1)Tr(g_2)$ disappears and this is the
standard central extension up to normalization.

{\bf Fact} \Hayashi88, \cite{GW89,FF92} 
The center of $U_{c=\kappa}(gl_n[t,t^{-1}],c)$
\footnote{more precisely one needs "local completion" of $U_{c=\kappa}(gl_n[t,t^{-1}],c)$, see \FFR}
is trivial: $\CC\times 1$,
unless $\kappa\ne -1$. For $\kappa=-1$ there exists large center. $c=-1$ is called "critical level".
%
%
{\Th \label{Z-th-app}  
\bea \label{center-fml}
:det^{col}(L_{full}(z)-\p) : ~~~\mbox{generates the center on the critical level}
\eea
i.e. define the elements $H_{ij}$ by
$:det^{col}(L_{full}(z)-\p) :=sum_{i=-\infty...\infty;~j=0...n} H_{ij} z^i \p^j$,
then $H_{ij}$ freely generate the center of $U_{c=-1}(gl_n[t,t^{-1}],c)$.
Here we consider $L_{full}(z)$ (see \ref{L-full}) as matrix valued generating function
for generators of $U_{c=-1}(gl_n[t,t^{-1}],c)$, i.e. $e_{mn}t^p\in U_{c=-1}(gl_n[t,t^{-1}],c)$.
}

The theorem was first proved in \CTBig, but  the definition of the normal ordering
given there is different. The coincidence with the standard normal ordering is due to
\cite{CM}, where another direct proof of the theorem will be given.
The proof in \CTBig~ is quite short, but requires various results: Talalaev's theorem, existence of the
center, AKS-arguments, remarkable ideas from \RybA,\RybB. 
In \cite{CM} direct proof will be given, as a corollary new proof of Talalaev's
theorem will be obtained.

\subsection{The geometric Langlands correspondence over
$\CC$}\label{GLCoC} The geometric Langlands correspondence over
$\CC$ \cite{BD} (excellent survey \Frenkel95) is related to the
constructions above as follows. (See \CTBig, which originates from \FFR,\GNR ~
  and is intimately related to Sklyanin's separation of variables).
As in number theory there are two versions: local and global.
\subsubsection{The  local correspondence:}
\bea
"Automorphic ~~side": & \mbox{Irreducible Representation }~~V~~ of~~ (gl_n[t,t^{-1}],c)_{c=crit} 
  \nn\\
 & \updownarrow Local~ Langlands ~   correspondence ~over ~\CC \updownarrow \nn\\
"Galois~~ side": &  \mbox{differential operator of order $n$:  } \sum_{k=0,...,n} \p^k h_k(z)
\eea
In number theory the local correspondence is a bijection between $Irreps(GL_n(\mathbb{Q}_p))$
and n-dimensional representations of the Galois group $Gal(\overline{\mathbb{Q}_p} , \mathbb{Q}_p)$,
where $\mathbb{Q}_p$ is the field of p-adic numbers.
The standard analogies: $\mathbb{Q}_p$ is analogous to $\CC((t))$, irreps of $GL_n(\mathbb{Q}_p)$ to
irreps of $(gl_n[t,t^{-1}],c)$, representations of the Galois group of $\mathbb{Q}_p$
to differentials operators (since  Galois group is similar to the fundamental group
(i.e. it is unramified part of the Galois group, when both are defined)
and a differential operators defines a representation of a fundamental group via monodromy of its solution).

Construction \CTBig:
Consider irrep $V$ of $(gl_n[t,t^{-1}],c)_{c=crit}$, by "Schur's lemma"
center of $U_{c=crit}(gl_n[t,t^{-1}],c)$ acts by scalar operators on any irrep, hence one
obtains a character $\chi_{V}$ of the center $\chi_V:ZU_{c=crit}(gl_n[t,t^{-1}],c)\to \CC$.
Now the desired differential operator is given by Talalaev's formula:
$\chi_V(:(det^{col}(\p-L_{full}(z)):)$. (See formula \ref{L-full} for $L_{full}(z)$ and \ref{center-fml}.)
Summary:
\bea
V\in Irreps(gl_n[t,t^{-1}],c) \leftrightarrows \chi_V: ZU_{c=crit}(gl_n[t,t^{-1}],c)\to \CC
\leftrightarrows \chi_V(:(det^{col}(\p-L_{full}(z)): ~ ). 
\eea
The non-triviality is hidden in Talalaev's formula and theorem \ref{Z-th-app}.
The correspondence is important in integrability since analysis of $(det^{col}(\p-L(z))Q(z)=0$
allows to write Bethe equations for the spectrum and defines eigenfunctions in separated variables.
It also satisfies many  geometric, categorical and "stringy" properties,
 but this is out of our scope.
\subsubsection{The  global correspondence }
Consider an algebraic curve $\Sigma$, denote by $ Bun_{GL_n}^{\Sigma}$ the moduli stack of $GL_n$ vector
bundles on $\Sigma$.
The original version of the correspondence \cite{BD} is the following
(later it was extended):
\bea
"Automorphic ~~side": & \mbox{Hitchin's D-modules on  } Bun_{GL_n}^{\Sigma}
  \nn\\
 & \updownarrow Global~ Langlands ~   correspondence ~over ~\CC \updownarrow \nn\\
"Galois~~ side": &  GL_n-opers ~ on ~\Sigma \nn \\
&  \mbox{ i.e. roughly speaking differential operators of order n on } \Sigma
\eea
The construction (not the proof) of the global case is actually  an immediate corollary of the local one (\cite{BD}).
It goes as follows. \\
Denote by $HH$ commutative algebra of differential
\footnote{actually twisted-differential operators} operators on $Bun_{GL_n}^{\Sigma}$
which are quantization of Hitchin's hamiltonians \Hitchin87.
Hitchin's D-modules $D$ are (by definition) parameterized by  characters $\chi_{D}:HH\to \CC$.
By nowadays standard arguments (i.e. using representation  $Bun_{GL_n}^{\Sigma}
= G_{in}\backslash GL_n((z))/G_{out}$)
one sees that the  commutative algebra $HH$ is actually a factor algebra of
$ZU_{c=crit}(gl_n[t,t^{-1}],c)$. So one gets $\tilde \chi_D:ZU_{c=crit}(gl_n[t,t^{-1}],c) \to \CC$
and by the local correspondence one gets the differential operator
$D_z$. One needs to check that it is actually globally defined object
on $\Sigma$ \cite{BD}. (It is quantization of a globally defined object - spectral curve,
so hopefully in future it  would be derived from  general facts).
\BX

{\Rem ~} For $\Sigma=\CC P^1$ with the "marked points" ("tame ramification") the Hitchin system (H.s.) is precisely
the standard Gaudin system with  Lax matrix \ref{Lax-st},
for marked points with multiplicities ("wild ramification") - \ref{Lax-glt}b
(see \cite{Markman}, \EnriquezRubtsovA95, \Nekrasov95). For $\Sigma=\CC P^1$ with nodes and cusps
see \CTHA, more general singularities \CTHB. H.s. can be seen as Knizhnnik-Zamolodchikov equation degenerated to the
critical level.
\subsubsection{Remarks}
Let us emphasize the two features of the integrability point of view on the correspondence.
\begin{itemize}
\item The correspondence has nothing specific to  Hitchin's integrable system,
it should work for an arbitrary integrable system which has a spectral curve.
The analog of $G$-oper should be "quantum spectral curve" defined by Talalaev's
type formula $det(\hat \lambda -L(\hat z))$. (For example in  Heisenberg  XXX-magnet
one should write $Det(1-e^{-\p} T(z))$, in super version: $Ber(1-e^{-\p} T(z))$, ...)
\item
The geometric Langlands correspondence over $\CC$ is just an appropriately understood quantization of the
basic correspondence in integrability:   Liouville tori of classical integrable system
correspond to some spectral curves. Quantization sends Liouville tori ($H_i=\l_i$) to Hitchin's
D-module (essentially one needs to construct $\hat H_i$, the D-modules are generated by $\hat H_i - \chi(\hat H_i)$,
for an arbitrary $\chi: \hat H_i\to \CC$)
and  (amusingly)  the spectral curve itself should
be quantized by Talalaev's type formula $det(\hat \lambda -L(\hat z))$.
The quantum spectral curve is a differential (difference - depending on system) operator in one variable.
The corresponding equation $det(\hat \lambda -L(\hat z)) Q(\hat z)=0$ generalizes
Baxter's T-Q relation. It is a key ingredient in a solution of quantum systems.
\end{itemize}
The geometric Langlands correspondence is rapidly developing in various directions,
much work should be done, one may hope for exciting achievements.
\subsection{
Capelli and Baxter type conjectures \label{Announce-s}} 
Let us announce some conjectures on  Talalaev's quantum spectral
curve, the Capelli identities and the Baxter type conjectures. We
plan to discuss this in subsequent papers.

\subsubsection{$det(\p-L(z))$ and  $det(\p+L(z))$ are conjugate  }
Let $L(z)$, $T(z)$ be  Lax matrices of the Gaudin and Yangian types
(see definitions \ref{GL-def}, page \pageref{GL-def} and  \ref{Yang-def}, page \pageref{Yang-def} respectively).
Under the formal conjugation $D\to D^*$ we mean the following: $(\sum_i f_i
(z)\p^i)^*= \sum_i (-\p)^i f_i (z)$ i.e. the conjugation does not
act on coefficients $f_i$.

{\em {\bf Proposition C1} The differential (difference) operators $det^{col}(\p-L(z))$
and  $(-1)^n det^{row} (\p+L(z))$
(respectively $det^{col} (1-T(z)e^{-\p})$ and $det^{row} (1-e^{\p}T(z))$)
are formally conjugate to each other. }

\SPRF The Yangian case follows from the known result, Gaudin case follows from the
previous one by degeneration. Indeed, obviously:
\bea
det^{col} (1-T(z)e^{-\p})= \sum_{k=0,...,n}   \sum_{j_1<...<j_k} qdet^{col}( T(z)_{j_1,....,j_k} ) (-e)^{-k\p} ,\\
det^{row} (1-e^{\p}T(z))= \sum_{k=0,...,n}   \sum_{j_1<...<j_k}  (-e)^{k\p} qdet^{row}( T(z)_{j_1,....,j_k} )   ,
\eea
where  $qdet^{col} (M(z))=\sum_{\sigma\in S_n} \prod_{k=1,...,n}
M_{\sigma(k),k}(z-k+1)$,
$qdet^{row} (M(z))=\sum_{\sigma\in S_n} \prod_{k=1,...,n}
M_{k, \sigma(k)}(z-n+k-1)$. For the Yangian type  matrix $M(z)$ it is true that:
$qdet^{col} (M(z))= qdet^{row} (M(z))$ (see e.g. formulas 2.24, 2.25 page 10 \Molev02).
Also recalling that any submatrix of the Yangian type matrix is again Yangian
type matrix we come to the desired conclusion for the Yangian case.

Now we follow the ideas of D. Talalaev \T04: one can represent $T(z) =1 + h L(z)+h^2(...)$,
hence: $1- e^{h \p }T(z)=  -h(\p + L(z) )+h^2(...) $. So $det^{col} (1-T(z)e^{-\p})=
(-h)^n det^{col} (\p+L(z)) +h^{n+1}(...)$.
Similar with $1- T(z) e^{ -h \p } $. We proved above that
$det^{col} (1-T(z)e^{-h \p})$ and $det^{row} (1-e^{h \p}T(z))$ are formally conjugate
(it was proved for $h=1$, but one can rescale with $h$ \T04).
So their first coefficient in expansion in $h$ are formally conjugate,
that is the desired result for the  Gaudin case. \BX

\subsubsection{Capelli type explicit formulas for $det(\p\pm
L(z))$  \label{Cap-fml-ss}} {\em {\bf Conjecture C2}  The Capelli
type explicit formulas. \label{Cap-expl-an-sect} Assume that
$\ms\times \ms$ Gaudin type Lax matrix $L(z)$ has only simple poles
at some  points $z=z_i,i=1,...,\mpn$ (e.g. formulas \ref{Lax-G-ex1},
\ref{Lax-st}), then: \bea
det^{col}(\p-L(z))= (det^{row}(\p+L(z)))^{formally~conjugate}= \nn \\ ~ \\ ~
~=
 \sum_{i=0,...,\ms} (-1)^{i} \Bigl( \sum_{
\begin{array}{c}
principal~i\times i \\ ~submatrices~ M~of~L(z)
\end{array}
}
det^{col} (M+ \sum_{v=1,...,\mpn} \frac{1}{z-z_v} diag \{
i-1,i-2,...,0 \} \Bigr) \p^{n-i} =\nn \\~~ \\~ = \sum_{i=0,...,\ms}
(-1)^{i} \p^{n-i} \Bigl( \sum_{
\begin{array}{c}
principal~i\times i \\ ~submatrices~ M~of~L(z)
\end{array}
 }
det^{col} (M+ \sum_{v=1,...,\mpn} \frac{1}{z-z_v} diag \{
0,-1,...,-i+1 \} \Bigr) = \nn \\~ \\~
%
\mbox{Now row-dets:
~~~~~~~~~~~~~~~~~~~~~~~~~~~~~~~~~~~~~~~~~~~~~~~~~~~~~~~~~~~~~~~~~~~~~~~~~~~~~~~~~~~~~~~~~~}
~ \nn \\~ =
 \sum_{i=0,...,\ms} (-1)^{i} \Bigl( \sum_{\begin{array}{c}
principal~i\times i \\ ~submatrices~ M~of~L(z)
\end{array}
}
det^{row} (M+  \sum_{v=1,...,\mpn} \frac{1}{z-z_i} diag \{
0,1,...,i-1 \} \Bigr) \p^{n-i} =\nn \\~~ \\~ = \sum_{i=0,...,\ms}
(-1)^{i} \p^{n-i}  \Bigl( \sum_{\begin{array}{c}
principal~i\times i \\ ~submatrices~ M~of~L(z)
\end{array}
}
det^{row} (M+  \sum_{v=1,...,\mpn} \frac{1}{z-z_i} diag \{
-i+1,...,-1,0 \} \Bigr). \nn \\ \eea

Analogously let us formulate conjecture on the explicit expressions
for the Nazarov type elements. Here determinants are substituted by
permanents.
\bea Tr( S^p (\p+L(z)))= \frac{1}{p!} \sum_{ 1 \le j_1 \le \ms, 1
\le j_2 \le \ms,..., 1 \le j_k \le \ms} \sum_{\sigma\in S_p}
\prod_{r=1,...,p} (L_{j_r, j_{\sigma(r)}}+ \sum_{v=1,...,\mpn}
\frac{(r-1)\delta_{j_r, j_{\sigma(r)}}}{z-z_v}  ) . \eea }

\subsubsection{General Capelli identities}
Let us denote by $L(z)$ Lax matrix of the Gaudin type given by the
formula \ref{Lax4Cap};  $T^{XXX}(z)$ is a Lax matrix given by
formula \ref{Lax-XXX-1}; recall that the linear map "Wick" was
defined by formula \ref{Wick-def}. Consider $H_{k,j}$ defined as
$det(\l-L(z))=\sum_{k=0,..n;j=-n,...,\infty} \l^k z^j H_{k,j}$. Any
function made of $H_{k,j}$ is called classical integral of motion
for the Gaudin system. Similar for the $T^{XXX}(z)$.

{\em {\bf Conjecture C3} For all classical integrals of motion $H_1,
H_2$ of the Gaudin and XXX systems: \bea
[ Wick (H^{classical~Gaudin}_1(p,q)), Wick (H^{classical~Gaudin}_2(p,q)) ]=0, \\
\nn\\
Wick( Tr (\l-L^{classical}(z))^k)=  Tr (\p-L^{quantum}(z))^k,\\
\nn\\
Wick ( Tr {S^k} (\l-L^{classical}(z)))= Tr {S^k} (\p-L^{quantum}(z)),\\
~ \nn\\ ~
\mbox{the Yangian type systems:} \nn \\
~
[ Wick (H^{classical~~XXX}_1(p,q)), Wick (H^{classical~XXX}_2(p,q)) ]=0, \\
~ \nn\\ ~
Wick( det (1-e^{-\l} T^{classical~XXX}(z)))= det (1-e^{-\p} T^{quantum~XXX}(z)),\\
~ \nn\\ ~
Wick( Tr (1-e^{-\l} T^{classical~XXX}(z))^k)= Tr (1-e^{-\p} T^{quantum~XXX}(z))^k,\\
~ \nn\\ ~
Wick ( Tr {S^k} (1-e^{-\l} T^{classical~XXX}(z)))= Tr {S^k} (1-e^{-\p}  T^{quantum~XXX}(z)),\\
~ \nn\\
~ \label{XXX-imman} \forall ~~Young~tableux:~~~ Wick (
(Pr_{Young~tableux} \one T^{classical}(z) \two T^{classical}(z) ....
)) = \nn \\ ~~~
 (Pr_{Young~tableux} \one T^{quantum~XXX}(z+c_1) \two T^{quantum~XXX}(z+c_2) ...  ). \\
\nn\\ ~\mbox{~~~~~~- No "Tr" generalizes \OkounkovB96, if one puts
"Tr" - weaker \Ok961} \nn \eea } {\Rem ~} The expressions
$(Pr_{Young~diagram} \one T^{quantum~XXX}(z+c_1) \two
T^{quantum~XXX}(z+c_2) ....  )$ are called "fused T-matrices" in
physics literature.
See \Ok961, \NazarovB96 ~
for the concise exposition and deep results for the particular case
of $U(gl_n)$, the elements in $ZU(gl_n)$ obtained in this way are
called "quantum immanants" in loc.cit. It seems there is no concise
mathematical exposition (and probably understanding) of general
"fusion".
$Pr_{Young~tableux}\in Mat_n\otimes ... \otimes Mat_n$ is the Young
projector on a irreducible $gl_n$ representation corresponding to a
Young tableux, numbers $c_i=row_i-column_i$ for the cell with a
content $i$.
The conjecture above generalizes the results from \Ok961, \OkounkovB96,  \MTVCap. 
%
\subsubsection{Generalized Baxter type conjectures}
{\em {\bf Conjecture C4} Let $(\pi,V)$ be a finite-dimensional
representation of $gl[t]$ or Yangian, then all  solutions of the
equations below are single-valued functions (products of rational
function over exponentials): \bea
 \pi(Tr (\p-L^{quantum}(z))^k) Q_1(z)=0,\\
\pi(Tr {S^k} (\p-L^{quantum}(z))) Q_2(z)=0 ,\\
\pi( Tr (1-e^{-\p} T^{Yangian}(z))^k) Q_3(z)=0,\\
\pi(Tr {S^k} (1-e^{-\p}  T^{Yangian}(z)))Q_4(z)=0, \eea and more
generally $\forall~ Young~tableux$: \bea \pi(Wick ( Tr
(Pr_{Young~tableux}  \stackrel{1}{\Bigl( 1-\l^{-1}
T^{classical}(z)\Bigr) }
\stackrel{2}{\Bigl( 1-\l^{-1} T^{classical}(z)\Bigr) } ....  ))) Q_5(z)=0,\\
\pi(Wick ( Tr (Pr_{Young~tableux}  \stackrel{1}{\Bigl( \l-
L^{classical}(z)\Bigr) }
\stackrel{2}{\Bigl(\l - L^{classical}(z)\Bigr) } ....  ))) Q_6(z)=0.
\eea }

This conjecture will hopefully lead to new forms of Bethe equations
for the spectrum, similar to the considerations with the standard
Baxter equations.

\end{document}